\documentclass[11pt]{article}
\usepackage{amsfonts,a4wide}
\usepackage{amssymb,soul}
\usepackage{graphicx}
\usepackage{mathtools}
\usepackage{enumerate}
\usepackage[dvipsnames]{xcolor}
\newtheorem{theorem}{Theorem}
 \newtheorem{corollary}[theorem]{Corollary}
 \newtheorem{lemma}[theorem]{Lemma}

 \newtheorem{proposition}[theorem]{Proposition}

\newtheorem{definition}[theorem]{Definition}
\newtheorem{remark}[theorem]{Remark}
\newtheorem{example}[theorem]{Example}

\newcommand{\RE}{\mbox{\rm Re}}
\newcommand{\diag}{\mbox{\rm diag}}

\newcommand {\proof} {\par{\it Proof}. \ignorespaces}
\newcommand {\eproof}
      {\space
        {\ \vbox{\hrule\hbox{\vrule height1.3ex\hskip0.8ex\vrule}\hrule}}
        \par}
  \DeclareMathOperator{\tr}{tr}
\DeclareMathOperator{\ii}{i}

\DeclareMathOperator{\p}{p}
\DeclareMathOperator{\dist}{dist}
\DeclareMathOperator{\Span}{span}
\DeclareMathOperator{\ran}{ran}
\DeclareMathOperator{\ap}{ap}
\DeclareMathOperator{\rev}{rev}

\newcommand{\Real}{\mathbb{R}}

\newcommand{\eps}{\varepsilon}

\newcommand{\To}{\Rightarrow}
\newcommand{\sap}{s_{\ap}}

\newcommand{\D}[1]{\mathcal{D}\left( #1\right)}

\newcommand{\set}[1]{\left\{#1\right\}}
\newcommand{\seq}[1]{\left<#1\right>}
\newcommand{\norm}[1]{\left\Vert#1\right\Vert}

\newcommand{\mycomment}[1]{}

\catcode`@=11
\font\tenex=cmex10 
\newdimen\p@renwd
\setbox0=\hbox{\tenex B} \p@renwd=\wd0 
\def\bmat#1{\begingroup \m@th
  \setbox\z@\vbox{\def\cr{\crcr\noalign{\kern2\p@\global\let\cr\endline}}%
    \ialign{$##$\hfil\kern2\p@\kern\p@renwd&\thinspace\hfil$##$\hfil
      &&\quad\hfil$##$\hfil\crcr
      \omit\strut\hfil\crcr\noalign{\kern-\baselineskip}%
      #1\crcr\omit\strut\cr}}%
  \setbox\tw@\vbox{\unvcopy\z@\global\setbox\@ne\lastbox}%
  \setbox\tw@\hbox{\unhbox\@ne\unskip\global\setbox\@ne\lastbox}%
  \setbox\tw@\hbox{$\kern\wd\@ne\kern-\p@renwd\left[\kern-\wd\@ne
    \global\setbox\@ne\vbox{\box\@ne\kern2\p@}%
    \vcenter{\kern-\ht\@ne\unvbox\z@\kern-\baselineskip}\,\right]$}%
  \null\;\vbox{\kern\ht\@ne\box\tw@}\endgroup}
\catcode`@=12
\newcommand {\mat}  [1] {\left[\begin{array}{#1}}
\newcommand {\rix}      {\end{array}\right]}

\catcode`@=11     
\font\tenex=cmex10 
\newdimen\p@renwd
\setbox0=\hbox{\tenex B} \p@renwd=\wd0 
\def\bmat#1{\begingroup \m@th
  \setbox\z@\vbox{\def\cr{\crcr\noalign{\kern2\p@\global\let\cr\endline}}%
    \ialign{$##$\hfil\kern2\p@\kern\p@renwd&\thinspace\hfil$##$\hfil
      &&\quad\hfil$##$\hfil\crcr
      \omit\strut\hfil\crcr\noalign{\kern-\baselineskip}%
      #1\crcr\omit\strut\cr}}%
  \setbox\tw@\vbox{\unvcopy\z@\global\setbox\@ne\lastbox}%
  \setbox\tw@\hbox{\unhbox\@ne\unskip\global\setbox\@ne\lastbox}%
  \setbox\tw@\hbox{$\kern\wd\@ne\kern-\p@renwd\left[\kern-\wd\@ne
    \global\setbox\@ne\vbox{\box\@ne\kern2\p@}%
    \vcenter{\kern-\ht\@ne\unvbox\z@\kern-\baselineskip}\,\right]$}%
  \null\;\vbox{\kern\ht\@ne\box\tw@}\endgroup}
\catcode`@=12    
\newif\ifMDlatex
\catcode`@=11     
\def\MD@us#1{\csname#1style\endcsname}
\def\MD@uf#1{\csname#1font\endcsname}
\def\MD@t{text}
\def\MD@s{script}
\def\MD@ss{scriptscript}
\newdimen\MD@unit
\MD@unit\p@
\def\MD@changestyle#1{
  \relax\MD@unit0.1\fontdimen6\MD@uf{#1}0
  \everymath\expandafter{\the\everymath\MD@us{#1}}
}
\def\MD@dot{$\m@th\ldotp$}
\def\MD@palette#1{\mathchoice{#1\MD@t}{#1\MD@t}{#1\MD@s}{#1\MD@ss}}
\def\MD@ddots#1{{\MD@changestyle{#1}%
  \mkern1mu\raise7\MD@unit\vbox{\kern7\MD@unit\hbox{\MD@dot}}%
  \mkern2mu\raise4\MD@unit\hbox{\MD@dot}%
  \mkern2mu\raise \MD@unit\hbox{\MD@dot}\mkern1mu}}%
\def\MD@iddots#1{{\MD@changestyle{#1}%
  \mkern1mu\raise \MD@unit\hbox{\MD@dot}%
  \mkern2mu\raise4\MD@unit\hbox{\MD@dot}%
  \mkern2mu\raise7\MD@unit\vbox{\kern7\MD@unit\hbox{\MD@dot}}}}%
\def\MD@vdots#1{\vbox{\MD@changestyle{#1}%
    \baselineskip4\MD@unit\lineskiplimit\z@
    \kern6\MD@unit\hbox{\MD@dot}\hbox{\MD@dot}\hbox{\MD@dot}}}%
\ifMDlatex
  \DeclareRobustCommand\ddots{\mathinner{\MD@palette\MD@ddots}}%
  \DeclareRobustCommand\iddots{\mathinner{\MD@palette\MD@iddots}}%
  \DeclareRobustCommand\vdots{\mathinner{\MD@palette\MD@vdots}}%
\else
  \def\ddots{\mathinner{\MD@palette\MD@ddots}}%
  \def\iddots{\mathinner{\MD@palette\MD@iddots}}%
  \def\vdots{\mathinner{\MD@palette\MD@vdots}}%
\fi
\catcode`@=12    

\usepackage{color} 

\begin{document}
\title
{Spectral theory of infinite dimensional dissipative Hamiltonian systems }
\author{C. Mehl \footnotemark[3]~\footnotemark[1]
\and V. Mehrmann\footnotemark[3]~\footnotemark[1]
\and  M. Wojtylak \footnotemark[2]~\footnotemark[4]
}
\maketitle

\begin{abstract}
The spectral theory for operator pencils and operator differential-algebraic equations is studied. Special focus is laid on singular operator pencils and three different concepts of singularity of operator pencils are introduced. The concepts are analyzed in detail and examples are presented that illustrate the subtle differences. It is investigated how these concepts are related to uniqueness of the underlying algebraic-differential operator equation, showing that, in general, classical results known from the finite dimensional case of matrix pencils and differential-algebraic equations do not prevail.  The results are then studied in the setting of structured operator pencils arising in dissipative differential-algebraic equations. Here, unlike to the general infinite-dimensional case,  the uniqueness of solutions to dissipative differential-algebraic operator equations is closely related to the singularity of the pencil. 
\end{abstract}

{\bf Keywords.} operator differential-algebraic equation, singular operator pencil, regular operator pencil,  dissipative Hamiltonian equation

{\bf AMS subject classification 2020.}
37l05,  37l20, 47d06,  93b28, 93c05
\noindent

\renewcommand{\thefootnote}{\fnsymbol{footnote}}
\footnotetext[3]{
Institut f\"ur Mathematik, MA 4-5, Technische Universit\"{a}t Berlin, Stra{\ss}e des 
17.~Juni 136,
D-10623 Berlin, Germany.
\texttt{$\{$mehl,mehrmann$\}$@math.tu-berlin.de}.
}

\footnotetext[2]{Instytut Matematyki, Wydzia\l{} Matematyki i Informatyki,
Uniwersytet Jagiello\'nski, Krak\'ow, ul. \L ojasiewicza 6, 30-348 Krak\'ow, Poland
   \texttt{michal.wojtylak@uj.edu.pl} (corresponding author).}
\renewcommand{\thefootnote}{\arabic{footnote}}

\section{Introduction}

\subsection{The setting} 
For regular and singular linear finite dimensional systems of differential-algebraic equations (DAEs) 
$
E\frac{d}{dt} x(t) =Ax(t)$, 
with real or complex $n\times m$ matrices $E,A$, the spectral theory is well established through the Kronecker 
canonical form of the matrix pencil $\lambda E-A$, see \cite{KunM06}. 
 In particular, the singular part of the Kroncecker canonical form (\cite{Gan59}) gives an exact description of the existence and uniqueness of solutions of the Cauchy problem $
E\frac{d}{dt} x(t) =Ax(t)$, $x(0)=x_0
$.

 In this paper we study the infinite-dimensional case and consider \emph{ 
linear operator 
differential-algebraic
equations  and the Cauchy problem,}
\begin{equation}\label{odae}
 E \frac {d}{dt} x(t)= A x(t),\quad x(0)=x_0\ \
\end{equation}
where the operator $E$  is a bounded operator from a Hilbert space $\mathcal X$ to a Hilbert space $\mathcal Y$,   $A: \mathcal D(A)\subset \mathcal X \to
\mathcal Y$ is densely defined and closed, and $x_0\in \mathcal X$. The variable $t$ will always stand for a real parameter. 
While in recent years several papers on this topic have appeared, see, e.g., \cite{CamM99,CamM97,EmmM13,ErbJMRT24,JacM22,LucSE99,PucRS18,Rei10,Tro18},
a systematic theory of singular operator pencils (i.e. infinite dimensional versions of Kronecker singular blocks) is still not well developed. 
In the present publication we fill this gap.  
In particular, we will be interested in the influence of the spectral properties of the operator pencil $\lambda E-A$  on the uniqueness of the solutions of \eqref{odae}.    
The systematic investigation of singular operator pencils leads to better understanding of general linear operator differential-algebraic equations, similarly as the theory of singular linear pencils is an important component of the theory of differential-algebraic equations and their  numerical solution, see \cite{KunM06}.  We show that this strategy is particularly useful in the case when additional structure is available in the operators $A,E$. 

The class of operator pencils that we consider is that arising from
the important class of  \emph{dissipative Hamiltonian differential-algebraic equations}, i.e., 
equations  of the form
\begin{equation}\label{dHODAEQ}
E\frac {d}{dt} x(t)= BQx(t),\ \ 
  x(0)=x_0,
\end{equation}
where $B$ is a dissipative operator   {and $Q^*E$ is selfadjoint and nonnegative}.

This  class of operator DAEs is arising in energy based modeling via port-Hamiltonian systems in almost all physical domains such as elasticity, electromagnetism, fluid dynamics, structural mechanics, geomechanics,
poroelasticity, gas or water transport, see e.g.
\cite{AltMU21,AltS17,AouCMA17,
BaaCEJLLM09,
BanSAZISW21,
CarMP17,
EggK18,
JacZ12,
Kot19,
KurZSB10,
LeGZM05,
MacM09,
MacMB05,
MacSM04a,
MacSM04b,
MatH13,
Mor24,
OetG97,
Ram19,
SchM02,
Vil07}.
An important feature of \eqref{dHODAEQ}  is the existence of an energy function (\emph{Hamiltonian})
$\mathcal H(x)$. This is often a quadratic function  and in the case of \eqref{dHODAEQ} given by
\[
\mathcal H(x)= \frac 12 \langle
Ex,Qx \rangle,
\]
where $\langle \cdot , \cdot \rangle$
is an inner product on $\mathcal Y$.
A special case of this equation is of the form
\begin{equation}\label{dHODAE}
E\frac {d}{dt} x(t)= (J-R)x(t),\ 
  x(t_0)=x_0,
\end{equation}
where $\mathcal Y=\mathcal X$, $B=J-R$, with $J=-J^*$ skew-adjoint,  $E=E^*$, $R=R^*$ self-adjoint and positive semidefinite, and $Q$ is the identity. 
The corresponding structured analysis has been a research topic of great interest,  see \cite{AchAM21,GerHR21,MehMW18,MehMW21,MehU23}, with surprising spectral properties compared to the general setting.
The analysis of these operator differential-algebraic equations is an active research area, see e.g.  the recent papers \cite{ErbJMRT24,GerT22,JacM22,MehZ23}. However, the main assumption in all these papers is the nonemptyness of the \emph{resolvent set}, i.e., the existence of $\lambda_0\in\mathbb C$ for which
$\lambda_0 E-A$ is boundedly invertible. The simple example  $A=E=\diag(1,1/2,1/3,\dots)$ in the space $\ell^2$ shows that this assumption is not necessary for existence and uniqueness of solutions of \eqref{dHODAE}, cf. Example~\ref{ex:37}.

In the following subsection we illustrate the need for a theory of singular operator pencils with a few examples.

\subsection{Examples}
We begin with a basic example which plays a crucial role in the first part of the paper. 

\begin{example}\label{ex:kblock}{\rm 
In the finite dimensional setting the canonical pencil having a right singular chain is a singular Kronecker block
\[
\lambda E-A= \mat{cccc} \lambda &-1  && \\ &\ddots&\ddots &\\ && \lambda&-1\rix .
\]
   {Observe that the matrix pencil
in this example is rectangular, and for any $\lambda_0\in\mathbb C$ the matrix $\lambda_0 E-A$ has an eigenvector. Hence, the resolvent set (set of regular points) is empty. It is also known that  the corresponding differential-algebraic equation $E\dot x = Ax$, $x(t_0)=x_0$ does not have a unique solution. In the paper we will take a closer look on the operator pencils with empty  set of regular points and corresponding differential-algebraic equations.  }

In the first part of the paper, when dealing with general unstructured operator pencils, generalizations of the above pencil will be studied that satisfy
\[
Ee_j=\alpha_j e_{j},\quad A e_j=\beta_j e_{j-1}.
\]
Here $e_j$ are canonical basis vectors of the space $\ell^2$ or $\ell^2(\mathbb Z)$ and $\alpha_j$, $\beta_j$ are complex scalars. Note that both operators $E$ and $A$ in this case may be unbounded, however they are both closed and densely defined. 
We will use this construction 
as a main source of counterexamples connected with operator differential-algebraic equations, see Examples \ref{ex:IA}, \ref{ex:noinf} and \ref{ex:36}.}
\end{example}

 While Example~\ref{ex:kblock} is a  toy-example, 
 we next present two realistic examples for  the  setting~\eqref{dHODAEQ}. 

\begin{example}\label{ex:stokes} {\rm
Equations of the form (\ref{dHODAE})  arise, e.g., in  fluid dynamics of incompressible Newtonian
fluids, see \cite{BatB00,ReiS23,Soh01,Tem77}. Let $v$ and $p$ denote the velocity and pressure, respectively,
considered as abstract functions, mapping the time $t\in \mathbb R$ into appropriate spatial function spaces,
\cite{EmmM13}. The leading order terms in the linearization of the instationary incompressible
\emph{Navier-Stokes equations}
\begin{align*}
\partial_t{v} - \nu \Delta v + (v\cdot\nabla) v + \nabla p &=
f \qquad \text{in } \Omega \times \mathbb R ,\\
\mbox{\rm~div~} v &=0 \qquad \text{in } \Omega \times \mathbb R,
\end{align*}
with constant density $\rho=1$ and viscosity $\nu=1$ in  the limit of the Reynolds number  going to $ 0$
leads to the instationary \emph{Stokes equation}
\begin{align}
\partial_t{v} - \Delta v +  \nabla p &= f
\qquad \text{in } \Omega \times \mathbb R,
\nonumber \\
\mbox{\rm~div~} v &=0 \qquad \text{in } \Omega \times \mathbb R.
\label{linns}
\end{align}
Omitting the functional analytic details, see e.g. \cite{EmmM13}, then formally
\begin{equation}\label{NSEJR}
E :(v,p) \mapsto (v,0),\ J:(v,p) \mapsto (-\nabla p,-\mbox{\rm~div~} v),\ R:(v,p) \mapsto (\Delta v,0),
\end{equation}
and it is clear that $E$ has as kernel the functions $(0,p)$ and that the pressure is only determined  up to
a constant.
The associated Hamiltonian in this example is $\frac 12 \langle v,v \rangle$ and does not depend on the pressure. 
The spectral theory of the associated operator and its use in the analysis and construction of numerical
methods is  currently an important research topic, see e.g. \cite{LauN21,ReiS23} and the references therein.
}
\end{example}

\begin{example}\label{ex:porousmedia}{\rm The analysis of linear poroelasticity in a bounded Lipschitz
domain $\Omega \subseteq \mathbb{R}^d$ with $d\in\{2,3\}$ (and boundary $\partial \Omega$) was introduced
in~\cite{Bio41}, see also \cite{Sho00}. In  \cite{AltMU21} a dissipative Hamiltonian mixed weak formulation
has been derived for the averaged displacement field $u$ and its time derivative $w$, as well as the
averaged pressure~$p$ that  satisfy 
{
\begin{align}
\label{eqn:twoField:opMatrix3}
\begin{bmatrix} Y & 0 & 0\\ 0 &  A_0& 0 \\ 0 & 0 &  M\end{bmatrix}
  \begin{bmatrix} \dot{w} \\\dot{u}\\ \dot{p} \end{bmatrix}
  = \begin{bmatrix} 0 & - A_0&  D^* \\  A_0& 0 & 0 \\ - D & 0 & - K
  \end{bmatrix}
  \begin{bmatrix} w\\ u\\ p \end{bmatrix}  + \begin{bmatrix} f\\ 0 \\ g \end{bmatrix}.
\end{align}}
In the original formulation the three operators~$ Y$, $ A_0$, and $ M$ on the
left-hand side are self-adjoint and positive definite and the Hamiltonian  associated with this system is
given by
\begin{equation}
	\label{eqn:energyfull}
	\mathcal{H}(w,u,p)
	\vcentcolon= \frac{1}{2}\, \Big(\langle  Y w, w\rangle + \langle  A_0u, u\rangle +
\langle M p, p\rangle\Big),
\end{equation}
where $\tfrac{1}{2}\langle  Y w,w\rangle$ describes the kinetic and $\tfrac{1}{2} \langle  A
u, u\rangle + \tfrac{1}{2}\langle  M p,p \rangle$ the potential energy.
When going over to the quasi-stationary solution one sets $ Y=0$. 
}
\end{example}

 We refer the 
 reader to \cite{MehZ23} for further applications in the model class
 \eqref{dHODAEQ}.

\subsection{An overview of the results}

The paper is organized as follows. In Section~\ref{sec:prelim} we introduce the notation and present some preliminary results concerning the finite and infinite spectrum of operator pencils.

Section~\ref{sec:sinpen} presents the general theory for singular operator pencils. While the situation in a finite dimensional space is clear, it appears that 
there is no natural analogue in infinite dimension.
In particular, the following conditions are equivalent in finite dimensions, while essentially different in the infinite dimensional situation:
\begin{enumerate}[\rm (a)]
    \item  the set of points for which $\lambda E-A$ is invertible is empty; 
    \item   {$\lambda E-A$} has a right or left singular polynomial: $(\lambda E-A)p(\lambda)=0$ or  $(\lambda E^*-A^*)p(\lambda)=0$  for some vector valued polynomial $p$; 
    \item   {$\lambda E-A$} has a right or left analytic holomorphic function: $(\lambda E-A)x(\lambda)=0$ or\\  $(\lambda E^*-A^*)x(\lambda)=0$  for some vector valued holomorphic function $x$;
    \item   {$\lambda E-A$} has a right or left approximate sequence of singular polynomials. 
\end{enumerate}
To deal with the situation we state five  natural conditions in Section~\ref{ss:wish} that a potential definition of singularity should satisfy: it should extend the finite dimensional definition, it should imply that there are no regular points, and it should be invariant under  taking reversal pencils and congruences, and finally, the orthogonal sum of two regular pencils should be regular. In subsequent Sections~\ref{sec:empty}--\ref{sec:appsing} we show that only (b) and (d)  satisfy these requirements.

In Section~\ref{sec:DAE} we study the relation between (a)--(d) and  the  (non)-uniqueness of solutions of the operator DAE \eqref{odae}.
Our main result in Theorem~\ref{nonunique} is the derivation of a sufficient condition for non-uniqueness of the solutions.
From this result it then follows that having a right singular chain implies non-uniqueness of the solutions, but not conversely. The remaining concepts (a), (c), and (d) seem to be unrelated to the question of (non)-uniqueness of solution in general, see examples in Section~\ref{sec:DAE}.

While the general situation seems initially unfavourable, it becomes dramatically different when considered in the  dissipative Hamiltonian \eqref{dHODAEQ} setting. In Section~\ref{sec:dhpencils} we give a complete characterization for the uniqueness of solutions,  see Theorems~\ref{th:ce} and~\eqref{tEJR}.
Several examples throughout the paper illustrate the results.

\section{Preliminaries}\label{sec:prelim}

\subsection{Basic notations}
In this paper we consider the spectral theory of differential-algebraic equations (DAEs) with coefficients that are operators in Hilbert spaces.
The scalar product and norm will be denoted respectively by $\seq{\cdot,\cdot}$ and $\norm{\cdot}$. It will be always clear from the context  in which space the norm is taken.  
 In several examples we will use the spaces $\ell^2=\ell^2(\mathbb N)$ and $\ell^2(\mathbb Z)$ of square-summable sequences. For the general theory of unbounded operators, only briefly presented below, we refer the reader e.g. to \cite{Wei12}.

By a linear operator (or, in short: operator) we understand a linear mapping $S:\mathcal D (S)\to
\mathcal Y$, where the \emph{domain} $\mathcal D(S)$ is a linear subspace of $\mathcal X$.
By $\mathbf B(\mathcal X,\mathcal Y)$ we denote the \emph{set of bounded linear operators} from $\mathcal X$ to $\mathcal{Y}$. If $S$ is a
closed, densely defined operator, then we say that it is   {\emph{boundedly invertible}}, if
there exists an operator $T\in{\mathbf B}(\mathcal Y,\mathcal X)$ such that $TS=I_{\mathcal D (S)}$ and $ST=I_{\mathcal
X}$, where $I_{\mathcal Z}$ denotes the identity operator on the space $\mathcal Z$.

 For   {a densely defined} operator $A:\mathcal D(A)\to \mathcal Y$ we define the  \emph{adjoint} operator $A^*$ as usually, i.e., we set $  {\D{A^*}}=\{y \in\mathcal Y:    \text{there exists } z\in \mathcal X\  \langle y, A x\rangle=\langle z , x \rangle ,\  x\in\mathcal D(A)\}$  and $  A^* y$ equals by definition the (unique) $z$ from the previous formula.  We will use without saying the fact that $A$ is closable if and only if $A^*$ is densely defined. In particular, in the space $\ell^2(\Gamma)$, where $\Gamma \in \{\mathbb N,\mathbb Z\}$, every operator for which the linear 
span of the standard orthogonal basis is contained in $\D A\cap\D{A^*}$
is closable and densely defined.

Although our focus will be on operator pencils associated with the operator DAEs (ODAEs) \eqref{odae} and~\eqref{dHODAE}, we will also deal with more general \emph{operator pencils} of the form
\begin{equation}\label{Pgeneral}
P(\lambda)=\lambda E-A,
\end{equation}
with  $E:\mathcal D(E)\to \mathcal Y$, $A:\mathcal D( A)\to \mathcal Y$,
  {$\mathcal D(E),\mathcal D(A)\subseteq\mathcal X$}, being closed and densely defined and such that  $P(\lambda)$ is closed and densely defined for all $\lambda \in\mathbb C$. This happens, in particular, if one of the operators $E$ or $A$ is bounded, however this is not the  only possibility.   We will always use the symbol  $\lambda$ for the free variable, therefore by $P(\lambda)$ will always mean the operator pencil in \eqref{Pgeneral}. When we  evaluate that pencil at a specific complex number $\lambda_0\in\mathbb C$ then $P(\lambda_0)$ denotes an operator from $\mathcal D( A-\lambda_0 E)$ to $\mathcal Y$.

When $\mathcal{X}_j$ ($j=1,2$) are Hilbert spaces, then by $\mathcal X_1\oplus\mathcal X_2$ we denote their orthogonal sum, i.e., their Cartesian product with the (unique) norm satisfying $\norm{ (f,g)}^2=\norm f^2+\norm g^2$. We will also use the symbol $f\oplus g$ for a pair $f\in\mathcal X_1$, $g\in\mathcal X_2$. If, additionally, $\mathcal{Y}_j$ ($j=1,2$) are Hilbert spaces and $S_j$  are operators from $\D{S_j}\subseteq\mathcal X_j$ to $\mathcal Y_j$ ($j=1,2$),  then by $S_1\oplus S_2$ we understand the operator from $\D {S_1}\oplus\D{S_2}$ to $\mathcal Y_1\oplus \mathcal Y_2$ given by $(S_1\oplus S_2)(f_1\oplus f_2)=S_1f_1\oplus S_2f_2$. In a similar way we denote the infinite orthogonal sum of Hilbert spaces $\bigoplus_j \mathcal X_j $, see e.g. \cite{Wei12} for the general theory.

To develop the spectral theory of operator pencils we first define various type of \emph{spectra}.

\pagebreak

\begin{definition}\label{def:spectra}
Consider an operator pencil as in \eqref{Pgeneral}.
\begin{enumerate}
\item By $s\big(P(\lambda)\big)$ we denote the \emph{set of singular points of $P(\lambda)$}, i.e., the set of all $\lambda_0\in{\mathbb C}$ for which the operator $ P(\lambda_0) $ is not   {boundedly invertible}. 
\item By $s_{\p}\big(P(\lambda)\big)$ we denote the \emph{set of point singularities of
$P(\lambda)$}, i.e., values $\lambda_0 \in \mathbb C$ for which there exists $x\in \mathcal D\big(P(\lambda_0)\big)\setminus\{0\}$ with $P(\lambda_0)x=0$.
\item By
$s_{\ap}\big(P(\lambda)\big)$ we denote the \emph{approximate singularities} of  $P(\lambda)$, i.e., the set of all points
$\lambda_0\in{\mathbb C}$ for which there exists a sequence $\{x_n\}_{n=1}^\infty$ of unit norm vectors 
$x_n\in \mathcal D\big(P(\lambda_0)\big)$ such that $P(\lambda_0)x_n\to0$. We will use without saying that the sequence $\{x_n\}_{n=1}^\infty$ above might be equivalently
assumed to satisfy $\,\liminf_{n\to\infty}\norm{x_n}>0$ only, instead of being normalized.
\item 
The set of \emph{regular points} is defined as  $\rho\big(P(\lambda)\big):={\mathbb C}\setminus
s\big(P(\lambda)\big)$.
\end{enumerate}
\end{definition}

These definitions of spectra of operator pencils are direct generalizations of respective definitions of spectra for operators when $\mathcal Y=\mathcal X$, $E=I_{\mathcal X}$. 
We highlight that all spectra defined in Definition~\ref{def:spectra} are viewed as subsets of the 
complex plane (excluding the point infinity). Infinity as a spectral point will be discussed in detail in Subsection~\ref{sec:infspec}.

\subsection{Finite spectrum of operator pencils} \label{sec:spect}

In this section we derive some properties of the spectra of general operator pencils,  extending some classical results for operators,  see e.g. \cite[Chapter 2]{Kub12}.
\begin{proposition}\label{prop:basic}
Let $P(\lambda)=\lambda E-A$ be as in \eqref{Pgeneral},
with $E\in\mathbf B(\mathcal X,\mathcal Y)$ and $A$ being a closed and densely defined operator from $\D A\subseteq\mathcal X$ to $\mathcal Y$.  Then the following statements hold:
\begin{enumerate}[\rm (i)]
\item  If $\lambda_0\in\rho\big(P(\lambda)\big)$ then the disc
\[
D(\lambda_0):=\set{z\in{\mathbb C}:|z-\lambda_0|< \frac1{\norm{(\lambda_0E-A)^{-1}}\cdot{\norm E}}}
\]
is contained in $\rho\big(P(\lambda)\big)$. In particular,
 $\rho(P\big(\lambda)\big)$  is an open subset of ${\mathbb C}$.
\item\label{l-boundary} If  $\rho\big(P(\lambda)\big)$ is nonempty, then the boundary of
    $s\big(P(\lambda)\big)$ 
is contained in $\sap\big(P(\lambda)\big)$.
\item\label{l-closed} The sets $s_{\ap}\big(P(\lambda)\big)$ and $s\big(P(\lambda)\big)$ are closed subsets of ${\mathbb C}$.
\end{enumerate}
\end{proposition}

\proof 
{  Note that for $E=0$ all statements become trivial (with the convention $1/0=\infty$), hence we assume that $E\neq 0$.  }

(i) For $\lambda_0 \in\rho\big(P(\lambda)\big)\cap{\mathbb C}$, the {   power series}
\[
\sum_{j=0}^\infty (-1)^j (z-\lambda_0)^j(\lambda_0E-A)^{-1}\big( E (\lambda_0E-A)^{-1}\big)^j,
\]
 {  convergent for $z\in D(\lambda_0)$, constitutes the inverse of the operator $P(z)$.}
Hence $D(\lambda_0)$ is contained in the set of regular points of of $P(\lambda)$.

(ii) Assume that $\lambda_0$ is on the boundary of $s\big(P(\lambda)\big)$ and let
 $\{\lambda_n\}_{n=1}^\infty\subseteq{\mathbb C}$ 
be a sequence contained in the set of regular points of $P(\lambda)$ that is converging to $\lambda_0$. 
The inclusion in (i) applied for each $n$ implies that
\begin{equation}
\dist\Big(\lambda_n,s\big(P(\lambda)\big)\Big)\geq\frac1{\norm{(\lambda_nE-A)^{-1}}\cdot{\norm E}}.
\end{equation}
As a consequence we obtain
\[
\norm{(\lambda_nE-A)^{-1}}\geq\frac1{\dist\big(\lambda_n,s\big(P(\lambda)\big)\big)\cdot\norm E}\geq
\frac1{|\lambda_n-\lambda_0|\cdot\norm E}\to \infty, \quad n\to\infty.
\]
Hence, there exists a sequence $\{f_n\}_{n=1}^\infty\subseteq\mathcal Y$ with $\norm{f_n}=1$ and
$\norm{(\lambda_nE-A)^{-1}f_n}\to\infty$. Then setting $g_n:=  (\lambda_nE-A)^{-1}f_n$
we have $g_n\in\mathcal D(\lambda_nE-A)=\mathcal D (A)$ and
\[
\norm{(\lambda_0 E-A)\frac{g_n}{\norm{ g_n}}}\leq |\lambda_0-\lambda_n|\norm E+\frac{\norm{(\lambda_n
E-A)g_n}}{\norm {g_n}}= |\lambda_0-\lambda_n|\norm E+\frac{\norm{f_n}}{\norm{g_n}}\to 0.
\]
Then, by definition, $\lambda_0\in s_{\ap}\big(P(\lambda)\big)$.

(iii) 
Let $\{\lambda_n\}_{n=1}^\infty\subseteq s_{\ap}\big(P(\lambda)\big)$ converge to some
$\lambda_0\in{\mathbb C}$.  Let also $\{x_n\}_{n=1}^\infty$ be a sequence of unit norm vectors in $\mathcal D
(A)$ such that $\norm{(\lambda_n E -A)x_n}\leq 1/n$ for $n\geq 1$. Such a sequence exists due to
$\lambda_n\in s_{\ap}\big(P(\lambda)\big)$. Since $E$ is bounded it follows that
\[
\norm{ (\lambda_0 E -A)x_n}\leq \norm{ (\lambda_n E -A)x_n} + |\lambda_n-\lambda_0|\cdot\norm{Ex_n}\to 0, \quad
n\to\infty,
\]
which shows that $\lambda_0\in s_{\ap}\big(P(\lambda)\big)$.
The closedness of $s\big(P(\lambda)\big)$ is clear by (i).
\eproof

\medskip

In this subsection we have reviewed several classical spectral concepts for linear operators and extended them to operator pencils. In the next subsection we discuss the spectrum at $\infty$. 
\subsection{Infinity as a spectral point, the reversal pencil}\label{sec:infspec}

In this section we use the concept of the reversal of an operator pencil to define and analyze the spectral properties at $\infty$.

\begin{definition} 
   \emph{The reversal  of an operator pencil $\lambda E-A$} is defined as the pencil $\lambda A-E$.  Further, we say that  $\infty$ is  a \emph{point singularity (approximate singularity, singularity, regular point)} of $P(\lambda)$ if zero belongs to the set of point singularities, (respectively approximate singularities,  singularities, regular points) of the reversal $\lambda A-E$. 
\end{definition}
 
It is clear that $\lambda_0\in\mathbb C\cup\{\infty\}$ is a singularity 
(approximate singularity, singularity, regular point) of 
$\lambda E-A$ if and only if $\lambda_0^{-1}$ is a singularity 
(approximate singularity, singularity,  regular point, respectively) of $\lambda A-E$, regardless whether the operators $A$ or $E$ are bounded or not, with the
standard conventions $\infty^{-1}=0$ and $0^{-1}=\infty$.
With this definition it is also tempting to consider the different singularities of the pencil as subsets of the extended complex plane $\mathbb C\cup\{\infty\}$. Unfortunately, however, Proposition~\ref{prop:basic} is not true if $\mathbb C$ is replaced by $\mathbb C\cup\{\infty\}$, as the following Example~\ref{ex:ce1} shows. This example also shows simultaneously  that Proposition~\ref{prop:basic} is not true if $A$ is assumed to be bounded instead of $E$.
\begin{example}\label{ex:ce1} {\rm
Let $A$ be a closed densely defined operator with spectrum equal to the whole complex plane, e.g., let $A$ be the multiplication operator by $z$ in $L^2(\mu)$, where $\mu$ is some finite measure supported on the whole complex plane. Thus, the operator pencil $\lambda I-A$ also has the spectrum  equal to the whole complex plane. The reversal  $Q(\lambda)=\lambda A-I$, however, has zero in the  set of its regular points as   {$Q(0)=-I$}. Therefore, infinity is a regular point for $P(\lambda)$ and the set of singular points is not closed in the extended complex plane $\mathbb C\cup\{\infty\}$. Similarly, the set of singular points of $Q(\lambda)$ equals $\mathbb C\setminus\{0\}$, which is not closed in $\mathbb C$.}
\end{example}

Therefore, in the following we keep viewing $s\big(P(\lambda)\big)$, 
$s_{\p}\big(P(\lambda)\big)$, $s_{\ap}\big(P(\lambda)\big)$, $\rho\big(P(\lambda)\big)$ 
as subsets of the complex plane, regardless on the behaviour of the pencil at infinity. 
A detailed study of infinity as a spectral point in the case that the resolvent set is nonempty is presented in the recent paper \cite{ErbJMRT24}.

\section{Singular operator pencils - general theory}\label{sec:sinpen}
The main aim of our paper is to study the situation when the set of regular points of an operator pencil
\eqref{Pgeneral} is empty, which can happen for several reasons.
Distinguishing between these different reasons requires a definition of a \emph{singular operator pencil}, see e.g. \cite{MehZ23} for some suggestions for possible definitions.  In order to see which definition is most appropriate  for a certain purpose, let us create a list of criteria for a good definition.

\subsection{Criteria for a definition of a singular operator pencil}\label{ss:wish}
We suggest that a definition of singularity for an operator pencil should meet
the following criteria.
\begin{enumerate}
\item In a finite dimensional situation it reduces to the standard definition. 
\item It implies that there are no regular points.
\item A pencil is singular if and only if its reversal is singular. 
\item If the pencils $\lambda E_j-A_j$ ($j=1,2)$ are not singular, then neither is $ (   \lambda E_1 -A_2 )\oplus ( \lambda E_2-A_2)$.

\item A transformation $S(\lambda E-A)T$, with $S,T$   {bounded and boundedly invertible}, keeps the pencil singular.
\end{enumerate}

  {Observe that item 3) implies that infinity is not a regular point. In particular, no pencil of the form $\lambda I-A$ should be called singular.
We believe that singularity of operator pencils should be a concept that can only be observed for operator pencils, but not for operators. This is in line with the 
finite dimensional case, where no pencil of the form $\lambda I-A$ can be singular.}

While the criteria 1)--5) seem quite   {straightforward}, we will see in the following that it is not easy to have them satisfied simultaneously.

\subsection{Emptyness of the set of regular points, initial remarks}\label{sec:empty}
Recall that for a general operator pencil $P(\lambda)$ as in \eqref{Pgeneral}  the set of point singularities $s_{\p}\big(P(\lambda)\big)$ is  
contained in the set of  approximate singularities $s_{\ap}\big(P(\lambda)\big)$, 
which in turn is contained in the set of all singularities $s\big(P(\lambda)\big)$, 
and we discuss all these sets as subsets of the complex plane, treating infinity 
separately. For this we consider the following three important cases:

\begin{eqnarray}
\label{s:spC}  s_{\p}\big(P(\lambda)\big)= \mathbb C  
&\text{ and }& \infty\text{ is a point singularity};\\
 \label{s:sapC}
s_{\ap}\big(P(\lambda)\big)= \mathbb C  &\text{ and }& \infty\text{ is an approximate singularity};
\\    \label{s:sC}
s\big(P(\lambda)\big)= \mathbb C & \text{ and }& \infty\text{ is a singularity}.
\end{eqnarray}
Clearly, \eqref{s:spC} implies \eqref{s:sapC} and this implies \eqref{s:sC},
and in the finite dimensional case it follows from the Kronecker canonical form \cite{Gan59} that all these conditions are equivalent.  We now present some examples that show that in the infinite dimensional case the
converse implications do not hold in general.

\begin{example}{\rm \label{IInotI}
Let $E$ be a bounded operator with $0$ in the approximate spectrum, but not in the point spectrum. Then
the pencil $P(\lambda)=\lambda E-E$ satisfies 
\eqref{s:sapC} but not \eqref{s:spC}, as, e.g., $P(2)=E$ has a trivial kernel.
}
\end{example}

\begin{example}{\rm\label{IIInotII}
Let $ A$ be a closed, densely defined operator on some Hilbert space $ \mathcal X$ having the complex plane
as its spectrum, but not as its approximate spectrum, e.g., $ A$ is a symmetric but not self-adjoint operator. Then
\begin{equation}\label{tildeA}
P(\lambda)=\lambda (I_{\mathcal X}\oplus 0) +  A\oplus 1,\quad \text{in } {\mathcal{ X}}\oplus {\mathbb C}
\end{equation}
is an operator pencil which satisfies  \eqref{s:sC} but not \eqref{s:sapC}.
}
\end{example}

These two examples are rather simple and do not present the nature of the problem yet.
More elaborate examples will be given subsequently.

  {Since in the finite-dimensional case singularity of a pencil is equivalent
to the condition that the spectrum equals $\mathbb C\cup\{\infty\}$},
one may wonder if any of the conditions \eqref{s:spC}--\eqref{s:sC} would be suitable as a definition of singularity for operator pencils, i.e., if this definition satisfies the requirements 1)--5)  of Subsection~\ref{ss:wish}.  The following example shows an immediate problems with item 4), see also \cite{KovP24} for a similar construction. 
\begin{example}\label{rem:notequal}\rm
It is well-known that the adjoint $S_0^\ast$ of the unilateral shift $S_0$ in $\mathcal X_0=\ell^2$ is a bounded linear operator
with the point spectrum being the open unit disc. Hence it is possible to find a
sequence of complex numbers
$\{\alpha_n\}_{n=1}^\infty$ such that
\[
\bigcup_{n=1}^\infty \sigma_p(S_0^\ast )+\alpha_n={\mathbb C}.
\]
Defining $\tilde {\mathcal X}=\bigoplus_{n=1}^\infty \mathcal X_0$ and $\tilde S=\bigoplus_{n=1}^\infty (
S_0^\ast +\alpha_n I_{\mathcal X_0})$, we get that the point spectrum of $\tilde S$ equals ${\mathbb C}$. Then the
pencil
\[
P(\lambda)=\lambda E-A \quad\text{with }E=I_{\tilde {\mathcal X}}\oplus 0,\; A:=\tilde S\oplus 1
\]
in $\mathcal X=\tilde {\mathcal{ X}}\oplus {\mathbb C}$
has the set of point singularities equal to $\mathbb C$ and infinity is a point singularity as well. 
Furthermore, it follows from the construction that one may represent a singular $P(\lambda)$ as an orthogonal  sum of two operator pencils, each of them having a nonempty set of regular points. As an example to see this, let
\[
J_1=\{n : \RE\,\alpha_n\leq 0\},\quad J_2=\{n:\RE\,\alpha_n> 0\}.   
\]
Setting ${\mathcal X_i}=\bigoplus_{j\in J_i}\mathcal X_0$  and $S_i=\bigoplus_{j\in J_i}(S_0^*+\alpha_j I_{\mathcal X_0})$, $i=1,2$, we obtain
\[
P(\lambda)=(\lambda I_{\mathcal X_1}- S_1)\oplus (\lambda I_{\mathcal X_2}-S_2)\oplus (\lambda 0-1),
\]
which (by combining one of the first two summands with the last) can obviously be reduced to two summands as required.
\end{example}

To obtain a definition of singular operator pencils satisfying our desired criteria, in the next subsection we generalize the concept 
of singular chains from the case of matrix pencils to operator pencils.

\subsection{Right and left singular polynomials}\label{sec:rightleft}
 In this subsection, we first review the theory of singular chains of matrix pencils.
Even in the finite-dimensional case, one has to distinguish between right and left singular chains and we will focus on the former ones first. 

We say that the $k+1$ vectors
$x_0,x_1,\dots,x_k$ ($k\geq 0$)  form a  right
singular chain for a
square complex matrix pencil $\lambda E-A$ if they  are linearly independent and satisfy
\begin{eqnarray}
\label{lchain}
&&Ax_0=0,\quad Ax_{j+1}=Ex_j\neq 0,\ j=0,\dots k-1,\quad Ex_k=0.
\end{eqnarray}
Note that for $k=0$  condition \eqref{lchain} just means that $\ker A\cap\ker E\neq\set0$.
Moreover, having a right singular chain directly implies \eqref{s:spC}. Indeed,  for any
$\lambda_0\in{\mathbb C}$ we can  define
\[
x(\lambda_0):=\sum_{j=0}^{k} \lambda_0^{j}x_{j}
\]
and,  due to the linear independence of $x_0,\dots, x_k$, we have that $x(\lambda_0)\neq 0$ and
\begin{eqnarray}\label{cal1}
(\lambda_0 E-A)x(\lambda_0)&=& \sum_{j=0}^{k}
\lambda_0^{j+1}Ex_{j}-\sum_{j=0}^{k}\lambda_0^{j}Ax_{j}\\
&=& \lambda_0^{k+1} Ex_k +\sum_{j=0}^{k-1}\lambda_0^{j+1}(Ex_{j}-Ax_{j+1}) + Ax_0=0.\label{cal2}
\end{eqnarray}
Furthermore, $\infty$ is a spectral point of $\lambda E-A$, because
$E$ has a nontrivial kernel.

Also let us recall, that if $x_0,\dots, x_k$ satisfy \eqref{lchain} but are not linearly independent or do not satisfy
$Ex_j\neq 0$ for all $j=0,\dots,k-1$, then a  shorter right singular chain is contained in the span of $x_0,\dots,x_k$.

In a similar way one defines left singular chains as right singular
chains of the conjugate transpose pencil $\lambda E^* -A^*$. Note that  every square matrix pencil that has
a left singular chain also  necessarily possesses  a right singular chain, although their
lengths may be different. 

Finally, a matrix pencil is called \emph{singular} if it possesses either a left or a right singular chain.  
From the discussion above, it  is apparent that one may equivalently define singular pencils as those for which there exists a nonzero vector valued polynomial $p(\lambda)$ with either $(\lambda_0 E-A)p(\lambda_0)=0$ for all $\lambda_0\in\mathbb C$ or $p^*(\lambda_0) (\lambda_0 E-A)=0$ for all $\lambda_0\in\mathbb C$. In view of these observations we introduce the following definition for operator pencils.

\begin{definition}\label{def:PS} \rm Consider an operator pencil of the form \eqref{Pgeneral}.
We say that an $\mathcal X$-valued polynomial $p(\lambda)$ is \emph{a right singular polynomial for the pencil $\lambda E-A$} if  $p(\lambda)$ is a nonzero polynomial, $p(\lambda_0)\in\D{E}\cap\D{A}$ and $(\lambda_0 E-A)p(\lambda_0)=0$ for every $\lambda_0\in\mathbb C$. 
We say that an $\mathcal Y$-valued polynomial $q(\lambda)$ is a left singular polynomial for the pencil $\lambda E-A$ if  it is a right singular polynomial for the adjoint pencil $\lambda E^*-A^*$. We call a pencil \emph{point singular} if has either a left or a right singular polynomial. 
\end{definition}

In the following we show that Definition~\ref{def:PS} satisfies the singularity criteria 1)--5) from Subsection~\ref{ss:wish}. 
For this we will need the following lemma   {on $\mathcal X$-valued polynomials. 
As in the standard case of scalar polynomials, we define the \emph{reversal}
of an $\mathcal X$-valued polynomial $p(\lambda)=\sum_{j=0}^k \lambda^j a_j$
($a_1,\dots,a_k\in\mathcal X$) to be the polynomial 
$\rev p(\lambda)=\sum_{j=0}^k \lambda^j a_{k-j}$ while a root of $p(\lambda)$ is
a value $\lambda_0\in\mathbb C$ satisfying $p(\lambda_0)=0$}.
\begin{lemma}\label{division}
Consider an operator pencil of the form \eqref{Pgeneral}.  If  the pencil $\lambda E-A$ has a right singular polynomial, then there exists a right singular polynomial $p_0(\lambda)$ such that 
both $p_0(\lambda)$ and $\rev p_0(\lambda)$ do not have any roots in the complex plane. Furthermore,  $\rev p_0(\lambda)$ is a right singular polynomial for $\lambda A-E$.
\end{lemma}
\proof
Suppose that $p(\lambda)=\sum_{j=0}^k \lambda^j a_j$   ($a_1,\dots,a_k\in\mathcal X$) is a right singular polynomial of $\lambda E-A$. As the linear span of $a_1,\dots,a_k$ is finite-dimensional, we can rewrite $p(\lambda)$ as  $p(\lambda)=\sum_{j=0}^l p_j(\lambda) f_j$, with some scalar polynomials $p_1(\lambda),\dots,p_l(\lambda)$, orthonormal vectors $f_1,\dots,f_l$ and $l\leq k$. 
Let $q(\lambda)$  be the greatest common divisor of $p_1(\lambda),\dots,p_l(\lambda)$. 
We set 
$$
p_0(\lambda)= \frac1{q(\lambda)}{ p(\lambda)}=\sum_{ j=1}^l  \frac{p_j(\lambda)}{q(\lambda)} f_j   .
$$
Observe that $p_0(\lambda)$ is an $\mathcal X$-valued polynomial, which has no zeros in the complex plane and is a right singular polynomial of $\lambda E-A$ due to $(A-\lambda E)\frac{p(\lambda)}{q(\lambda)}=\frac 1{q(\lambda)} (A-\lambda E){p(\lambda)}=0$. By the properties of the reversal polynomial, $\rev p_0(\lambda)$ has  no zeros in the complex plane as well and is a right singular polynomial for the reversal pencil $\lambda A-E$.
\eproof

\medskip

We now show that a point singular operator pencil has no regular points.
\begin{proposition}\label{polyC} Consider an operator pencil $P(\lambda)=\lambda E-A$ 
with $E\in\mathbf B(\mathcal X,\mathcal Y)$ and $A$ being a closed and densely defined operator from $\D A\subseteq\mathcal X$ to $\mathcal Y$.  Then the following statements hold:
\begin{enumerate}[\rm (i)]
    \item 
If $\lambda E-A$ has a right singular polynomial $p(\lambda)$, then the set of point singularities coincides with the complex plane and also $\infty$  is a point singularity, i.e., \eqref{s:spC} holds.
\item
If $\lambda E - A$  has a left singular polynomial,  then the set of all singularities of $\lambda E-A $ is the whole complex plane and also  $\infty$ is a singularity, i.e., \eqref{s:sC} holds.
\end{enumerate}
\end{proposition}
\proof
 (i) The first statement is obvious by Lemma~\ref{division}. 
Part (ii)   follows directly from the fact that if an operator  $Z^*$ is not   {boundedly invertible}, then neither is $Z$.
\eproof

\medskip

We now have the following observations concerning the singularity criteria from Subsection~\ref{ss:wish}. The concept of point singularity in Definition~\ref{def:PS}  obviously satisfy criteria 1) and 5). Criterion 2) follows from Proposition~\ref{polyC}, and criterion 3) from Lemma~\ref{division}.  Criterion 4) 
follow by contraposition from the fact that if
\[
\big((  \lambda E_1-A_1)\oplus ( \lambda E_2-A_2)\big) \ \big(p_1(\lambda)\oplus p_2(\lambda)\big)  =0,
\]
then $p_j(\lambda)$ is either the zero polynomial or a right singular polynomial for $\lambda E_j-A_j$ $(j=1,2)$.
%
\begin{example}\label{ex:nodecom}{\rm
The operator pencil $\lambda E-A$ of Example~\ref{rem:notequal} has $\mathbb C$ as set of point singularities and infinity is a point singularity as well.  However, the operator pencil is  not point singular according to Definition~\ref{def:PS}. Indeed, from the construction it follows that it can be decomposed into two pencils with nonempty  sets of regular points.}
\end{example}

Although having a right or left singular polynomial is in accordance with our list of singularity criteria, 
in the infinite dimensional context it is rather restrictive. 
In the following subsections we will therefore discuss other  concepts.

\subsection{Right and left singular holomorphic functions}\label{ss:singularanalytic}

{  Note that a right singular polynomial evaluated at  $\lambda_0$ is an 
eigenvector of $\lambda_0 E -A$. An obvious extension of this concept is obtained via the
transition to holomorphic functions, defined except for some small set, so that it provides a  set of point singularities dense in $\mathbb C$. Surprisingly, we will observe that such an extension of the definition appears to be  not suitable as a concept describing singularity of an operator pencil.}

\begin{definition}\rm
Consider an operator pencil of the form \eqref{Pgeneral}. We say that an 
$\mathcal X$-valued function $x(\lambda)$ is a \emph{left singular function} for
$\lambda E-A$ if it is defined and holomorphic on $\mathbb C$ except, possibly, for a 
discrete set of points and 
\[
(\lambda E- A)x(\lambda)= 0, \quad x(\lambda)\neq 0,
\]
for all $\lambda\in\mathbb C$,  except, possibly, for a discrete set of points. 
\end{definition}

For completeness we state an analogue of Proposition~\ref{polyC}, noting that for the proof of (i) one needs to use the fact that the set of zeros of a  holomorphic function  is nowhere dense in $\mathbb C$.
\begin{proposition} Let $P(\lambda)=\lambda E-A$,
with $E\in\mathbf B(\mathcal X,\mathcal Y)$ and $A$ being a closed and densely defined operator from $\D A\subseteq\mathcal X$ to $\mathcal Y$.  Then the following statements hold:
\begin{enumerate}[\rm (i)]
    \item 
If $\lambda E-A$ has a right singular function $x(\lambda)$, then the set of point singularities contains   all nonzero points of the domain of  the function $x(\lambda)$. In particular, the set of approximate singularities is the whole complex plane.
\item
If $\lambda E - A$  has a left singular function,  then the set of all singularities of $\lambda E-A $ is the whole complex plane.
\end{enumerate}
\end{proposition}

We now present two examples of pencils having right singular functions. 

\begin{example}\label{ex:IA}{\rm
Consider  the following operators in the Hilbert space $\ell_2$ of 
square-summable sequences:
\[
\lambda E-A = \mat{ccccc} \lambda & -1  &&&\\ &\lambda/2 &-1 &&\\ && \lambda/3 & & \\ &&& \ddots &\ddots \rix,
\]
i.e., $A$ is the backward shift, and $E$ is a diagonal operator. With $x_j=\frac1{(j-1)!}e_j$, $j=1,2,\dots$, where $e_1,e_2,\dots,$  stands for the canonical basis of $\ell^2$, one clearly has that 
\[
x(\lambda)=\sum_{j=1}^\infty \lambda^j x_j
\]
is everywhere square-summable and nonzero. Furthermore, 
\[
Ax_1=0\quad\mbox{and}\quad Ax_j=\frac 1{(j-1)!} e_{j-1}=Ex_{j-1},\quad j=2,3,\dots 
\]
Therefore, $\sum_{j=1}^\infty \lambda^j A x_j$ is square-summable and $x(\lambda)$ belongs to $\D A$ for all $\lambda$. Moreover, we have 
\[
(\lambda E-A)x(\lambda)= \lambda \sum_{j=1}^\infty \lambda^j Ex_j
-\sum_{j=1}^\infty \lambda^j Ax_j=0,
\]
i.e. $x(\lambda)$ is a right singular function of $\lambda E-A$.
}
\end{example}

\begin{example}\label{ex:noinf} {\rm
Consider a pencil $\lambda I -A$ in the space $\ell^2(\mathbb Z)$, i.e., the
Hilbert space of square-summable sequences with index set $\mathbb Z$, where the operator $A$ is defined with the use of the canonical basis as follows:
\[
A_0(e_j)= \frac{|j|!}{|j-1|!}e_{j-1},\quad j\in\mathbb Z.
\]
Clearly, $A_0$ is a closable, densely defined operator, and let $A$ be its closure. Define $x(\lambda)$ by the Laurent series
\[
x(\lambda)= \sum_{j=-\infty}^\infty \frac{\lambda^j}{|j|!} e_j,
\]
which is convergent everywhere except at zero and infinity. Furthermore, 
\[
A   \frac{e_j}{|j|!} =\frac 1{|j-1|!} e_{j-1},\quad j\in\mathbb Z,
\]
and therefore $\sum_{j=-\infty}^\infty \lambda^j A  \frac{e_j}{|j|!}$ is square-summable and $x(\lambda)$ belongs to $\D A$ for all $\lambda\neq 0$.  Moreover, $(\lambda I-A)x(\lambda)= 0$ as it can be immediately seen by shifting the summation index in $Ax(\lambda)$.
}
\end{example}

In view of Example~\ref{ex:noinf},  `having a right or left singular function' does not seem to be a good  generalization of the notion of a singular matrix pencil. It satisfies the criteria 1), 2), 4) and 5) from Subsection~\ref{ss:wish}, 
but it does not satisfy criterion 3). One could think of some variations, e.g., `both $\lambda E-A$ and $\lambda A-E$ have a (left or right) singular function' as another definition of singularity. This would satisfy the criteria 1), 2), 3) and 5), but no longer criterion 4) as the following example demonstrates.

\begin{example} \label{ex:non4}{\rm Let $A$ be the operator from Example \ref{ex:noinf} and consider the operator pencil
\[
\lambda (I\oplus A) - (A\oplus I)= (\lambda I-A)\oplus(\lambda A-I). 
\]
Then $x(\lambda)\oplus 0$ is a right singular function for the pencil, while $0\oplus x(\lambda)$ is a right singular function for the reversal.
It should be noted, that none of the coefficients of the pencil is bounded, however, the pencil is closed and densely defined for all $\lambda\in\mathbb C$.
}
\end{example}
\subsection{Operator pencils with approximate joint kernel}\label{sec:jointkernel}
One of the situations that we would like to cover by the notion of singularity is the case of operator pencils of the form~\eqref{Pgeneral} where $A$ and $E$
have an \emph{approximate joint kernel}, i.e., 
there exists a sequence of unit norm vectors $\{x_n\}_{n=1}^\infty\subseteq\D A\cap\D E$ with $Ax_n\to 0$
and $Ex_n\to 0$ for $n\to \infty$. Clearly, if this happens then $s_{\ap}\big(P(\lambda)\big)=\mathbb C$ and infinity is an approximate singularity as well. {  Further, it is clear that for any $\eps>0$ there exist $A_1$, $E_1$, both of norm bounded by $\eps$, such that $\lambda(E+E_1)-(A+A_1)$ is a point singular pencil.}

Therefore, it is tempting to call such pencils singular. Note that this situation is not covered by the concept of point singularity in Definition~\ref{def:PS}, unless the dimension is finite, as the following example demonstrates. 
\begin{example}\label{ex:diag}{\rm
    Let $\mathcal X=\mathcal Y=\ell^2$ and let $A=E=\diag(1,1/2,1/3,\dots)$. Here both operators $E$ and $A$ are bounded (even compact) and the pencil $\lambda E-E$ has a unit norm sequence $x_n:=e_n$  (the $n$-th vector of the canonical basis) with $Ax_n=Ex_n\to 0$. }
\end{example}

\begin{remark}\label{rem:debate}{\rm
In some publications there is a debate whether an operator pencil as in  Example~\ref{ex:diag} should be considered as singular. This results from the fact that the corresponding linear relation
\begin{equation}\label{relation}
    \{ (Ex,Ax) : x\in\mathcal X\}\subseteq \mathcal X\times \mathcal X
\end{equation}
is the identity relation.
It was shown in \cite{BerSTW22,BerTW16} that in finite dimensional spaces there is a one-to-one correspondence between linear relations and linear pencils that keeps the spectral properties. In particular, if $\mathcal X=\mathcal Y$ is finite dimensional, then any linear pencil $\lambda E-A$ with the corresponding relation \eqref{relation} being equal to the identity is a regular pencil, congruent to $\lambda I-I$. 
Apparently, in the infinite dimensional setting the correspondence is  more complex.
See \cite{GerT22,KhlTW23} for further work on regular and singular linear relations in infinite dimensional spaces.
    }
\end{remark}

{  We now discuss the relation between approximate joint kernels and the distance to singularity of finite sections of operator pencils. 
The latter are defined as follows.} 
Consider two sequences of (not necessarily
orthogonal) projections  $P_n:\mathcal X\to \mathcal X_n$, $Q_n:\mathcal Y\to \mathcal Y_n$ onto finite dimensional subspaces $\mathcal
X_n\subseteq  \D A\cap \D E$, $\mathcal Y_n\subseteq\mathcal{Y} $ 
  {of dimension $k_n$} that satisfy the following conditions:
\begin{eqnarray}
&&P_n f\to f\ \mbox{for }n\to\infty,\quad f\in \mathcal X,  \label{P1}\\
&&Q_n A P_n f\to Af,\ Q_n E P_n f\to Ef \ \mbox{for }(n\to\infty),\   f\in\D A\cap \D E.  \label{P3}
\end{eqnarray}
This is the typical setting in numerical methods where the solution of
partial-differential equation is approximated using some discretization method, e.g. via Petrov-Galerkin projection
in the finite element method, or truncated Taylor  or Fourier series, 
see e.g.
\cite{BofBF13,Boe18,BoeM20,Cia02} and the references therein.

{   While  the operators $Q_nEP_n$ and $Q_nAP_n$ act from
$\mathcal X$ to $\mathcal Y$, we will also consider
$Q_nE\vert_{\ran P_n}$ and $Q_nA\vert_{\ran P_n}$
acting from $\mathcal X_n$ to $\mathcal Y_n$, both being $k_n$-dimensional. Choosing orthonormal bases in these subspaces (with respect to original inner products on {  $\mathcal X$} and $\mathcal Y$), we can identify the two operators with $k_n\times k_n$ matrices and use the Frobenius norm
$\norm{X}_F:=(\tr (X^*X))^{1/2}$ on $\mathbb C^{k_n\times k_n}$. Thus, for fixed $n$ we can measure the \emph{distance to singularity} of $Q_n(\lambda E-A)\vert_{\ran P_n}$  as}

%
\begin{equation}\label{dsin}
  {\delta_n}:=\inf\set{{  \sqrt{\big\|Q_n E\vert_{\ran P_n} -\tilde E\big\|^2_F +\big\|Q_n A\vert_{\ran P_n}-\tilde A\big\|^2_F}:  {\tilde E,\tilde A\in\mathbb C^{k_n\times k_n},}\;\lambda \tilde E-\tilde A\text{
singular }}},
\end{equation}
as in \cite{ByeHM98}.
\begin{proposition}\label{distto0}
Consider a pencil $\lambda E-A$ of the form \eqref{Pgeneral}. If there exists a sequence of unit norm
vectors $\{x_n\}_{n=1}^\infty$ contained in the domain of $A$ and such that $Ax_n\to 0$, $Ex_n\to 0$ for $n\to\infty$,
then for any sequence of projections $P_n,Q_n$ satisfying \eqref{P1}--\eqref{P3}, the distance to singularity \eqref{dsin}
of the matrix pencils
$ Q_n ( \lambda E- A  )\vert_{\ran P_n}$ converges to zero with $n\to\infty$.
\end{proposition}
\proof
Let $\eps>0$ be arbitrarily small and let $n$ be such that $\norm{Ax_n}<\eps$, $\norm{Ex_n}<\eps$. Moreover,
let $m_0$ be such that for all $m>m_0$ one has
$\norm{P_mx_n}>1-\eps$, $\norm{Q_mEP_mx_n}<2\eps$, $\norm{Q_mAP_mx_n}<2\eps$. The existence of such an $m_0$ 
follows by  \eqref{P1} and \eqref{P3}.
Then with $g_m=\frac{P_mx_n}{\norm{P_mx_n}}\in\ran P_m$, $\norm g_m=1$ for $m>m_0$, one has
\[
\sigma_{\min}^2\left(\mat{c}Q_m A\vert_{\ran P_m}\\ Q_m E\vert_{\ran P_m}\rix\right)\leq \norm{ Q_m A g_m  }^2+\norm{Q_mE
g_m}^2<\frac{8\eps^2}{(1-\eps)^2}.
\]
As $\eps>0$ was arbitrary, the left hand side has to converge to zero with $m\to \infty$, which by an
estimate in  \cite{ByeHM98} implies that \eqref{dsin} converges to 0 with $n\to\infty$.
\eproof
\medskip
{  We remark that the  convergence of \eqref{distto0} for (some)  $P_n$, $Q_n$ satisfying   \eqref{P1}--\eqref{P3} may, nevertheless, happen even for regular pencils.}
\begin{example}\label{ex13}{\rm
Let $T$ be the bilateral shift in $\ell^2(\mathbb Z)$, i.e., $T$ acts on the canonical basis as
$Te_j=e_{j+1}$ for $j\in\mathbb Z$. As $T$ is unitary, for the pencil $P(\lambda)=\lambda T-T$ we have
$s\big(P(\lambda)\big)=\{1\}$ and \eqref{s:sapC} does not hold and so the pencil should not be seen as singular in any sense. However, the projections $P_n=Q_n$ onto
$\Span\{  {e_{1-n}},\dots,e_0,e_1,\dots,e_n\}$ satisfy \eqref{P1}--\eqref{P3} and the pencil
$
P_n (\lambda T-T) P_n
$
is clearly singular as $e_n\in\ker P_nTP_n$ for $n\geq  1$. }
\end{example}

\subsection{Approximate singularity}\label{sec:appsing}
In this subsection we define and analyze the concept of approximate singularity for operator pencils.

\begin{definition}\label{def:apppol}{\rm
A sequence of $\mathcal X$-valued polynomials $\{ p_n(\lambda)\}_{n=1}^\infty$ is called a \emph{right approximate polynomial sequence   for the operator pencil $\lambda E-A$} of the form \eqref{Pgeneral} if 
\begin{equation}\label{pnot0}
p_n(\lambda_0)\neq 0,\quad
     \lim_{n\to \infty}\|p_n(\lambda_0)\|\neq 0,\quad\rev p_n(\lambda_0)\neq 0,\quad  \lim_{n\to\infty}\|\rev p_n(\lambda_0)\|\neq 0,
\end{equation}
while $p_n(\lambda_0)\in\D{E}\cap\D{A}$  for $n\geq 1$ and 
\[
\lim_{n\to \infty}\norm{(\lambda_0 E-A)p_n(\lambda_0)}=0, \  \lim_{n\to \infty}\norm{(\lambda_0 A-E)\rev p_n(\lambda_0)}=0,
\]
for all $\lambda_0\in\mathbb C$.
We say that a sequence  $\{ q_n(\lambda)\}_{n=1}^\infty$ of $\mathcal Y$-valued polynomials   is a \emph{left  approximate polynomial sequence   for the pencil $\lambda E-A$} if  it is a  right approximate polynomial sequence   for the  adjoint pencil $\lambda E^*-A^*$. Finally, we say that a pencil is \emph{approximately singular} if it has either a left or right approximate polynomial sequence. }
\end{definition}
One may wonder if $p_n(\lambda_0)\not\to 0$ for $\lambda_0\in\mathbb C$ implies $\rev p_n(\lambda_0)\not\to 0$ for $\lambda_0\in\mathbb C$, or, in other words, if the conditions on the reversal in Definition~\ref{def:apppol} are necessary. The following example demonstrates the necessity.
\begin{example}\rm Let
$p_n(\lambda)= e_1+ (\lambda^n/n!) e_2$ be a sequence of $\mathbb C^2$-valued polynomials. Clearly $p_n(\lambda_0)$ does not converge to 0 for any $\lambda_0\in\mathbb C$, while $\rev p_n(\lambda_0)=\lambda^n e_1 + e_2/n!$ does converge to 0 for $|\lambda_0|<1$. 
\end{example}

It is obvious that the situation of having a joint approximate kernel as in Subsection~\ref{sec:jointkernel} is covered by the definition of approximate singularity, but (even in the finite dimensional case) the converse is not true. We have the following lemma.
\begin{lemma}\label{xns} Consider an operator pencil $P(\lambda) =\lambda E-A$ of the form \eqref{Pgeneral}.
\begin{enumerate}[\rm (i)]
   \item If $P(\lambda)$ has a right approximate polynomial sequence, then there exist sequences of unit norm vectors $\{x_n\}_n\subseteq \D A$ with $Ax_n\to 0$ and $\{y_n\}_n\subseteq \D E$, with $Ey_n\to 0$, for $n\to\infty$.
 \item If $P(\lambda)$ has a sequence of unit norm vectors $\{x_n\}_n$ for which $Ax_n\to 0$ and $Ex_n\to 0$, then it has a right approximate polynomial sequence.
   \end{enumerate}
\end{lemma}
\proof For the proof of (i) set $x_n=\frac{p_n(0)}{\|p_n(0)\|}$, 
$y_n=\frac{\rev p_n(0)}{\|\rev p_n(0)\|}$, $n>1$. For the proof of (ii) take constant polynomials $p_n(\lambda)=x_n=\rev p_n(\lambda)$.
\eproof

\medskip

In Section~\ref{sec:dhpencils}  we present the class of dissipative Hamiltonian operator pencils for which the condition in Lemma~\ref{xns}(ii) is already equivalent  to being approximately singular.

Next we show that, among other properties,  for the concept of approximate singularity the criteria 1)--5) from
Subsection~\ref{ss:wish} are satisfied. Let us start with basic properties, that contain 1).
\begin{proposition}
Consider an operator polynomial  $P(\lambda)=\lambda E-A$ of the form \eqref{Pgeneral},
with $E\in\mathbf B(\mathcal X,\mathcal Y)$ and $A$ being closed and densely defined operator from $\D A\subseteq\mathcal X$ to $\mathcal Y$.  Then the following statements hold:
\begin{enumerate}[\rm (i)]
   \item  If $P(\lambda)$ has a    {right singular polynomial} then it has a right approximate polynomial sequence.
   \item  If the spaces $\mathcal X$ and $\mathcal Y$ are finite dimensional then a pencil is approximately singular if and only if it is singular. 
   \end{enumerate}
\end{proposition}

\proof Item (i) follows directly from Lemma~\ref{division}. To see the statement in  (ii) observe that if $\lambda E-A$ has an approximate polynomial sequence, it is not   {boundedly invertible} for any $\lambda\in\mathbb C$. Hence, if $\mathcal X$ and $\mathcal Y$ are finite dimensional, the pencil is singular.
\eproof

\medskip

Let us now move to criterion 2) from Subsection~\ref{ss:wish}. We state the following proposition for completeness of the presentation, the proof is obvious.

\begin{proposition}\label{PP?}
Consider an operator pencil $P(\lambda)=\lambda E-A$ of the form \eqref{Pgeneral},
with $E\in\mathbf B(\mathcal X,\mathcal Y)$ and $A$ being a closed and densely defined operator from $\D A\subseteq\mathcal X$ to $\mathcal Y$.  Then the following statements hold:
\begin{enumerate}[\rm (i)]
    \item  If $\lambda E-A$ has right approximate polynomial sequence, then the set of approximate singularities coincides with the complex plane and infinity is an approximate singularity.
\item
If $\lambda E-A$ has left approximate polynomial sequence,  then the set of all singularities of $\lambda E-A $ coincides with the whole complex plane and infinity is an approximate singularity.
\end{enumerate} 
\end{proposition}

We discuss now the remaining criteria of Subsection~\ref{ss:wish}. Criterion 3) follows directly from the definition. 
Criterion 4) follows from the fact that  $x_n\oplus y_n$ converges to zero if and only if both $x_n$ and $y_n$ converge to zero. 
Criterion 5) is obvious. 
\begin{example}\label{ex:nondecom2}{\rm
The operator pencil $\lambda E-A$ of Example~\ref{rem:notequal} has $\mathbb C$ as set of approximate singularities and infinity is an approximate singularity as well.  (In fact,
these are already point singularities.) However, the pencil is  not approximately singular according to Definition~\ref{def:apppol}. Indeed, this would contradict Criterion 4), as the pencil can be decomposed into two pencils with nonempty  sets of regular points.}
\end{example}

In the following we present further examples. In order to check 
condition~\eqref{pnot0}, a useful tool is provided by the following proposition. 
\begin{proposition}\label{PP}
For $p_n(\lambda)=\sum_{j=0}^{k_n} \lambda^j x^{(n)}_j$ with $x^{(n)}_j\in\mathcal X$, let     
$$\xi:=\inf_{n\geq 0}\lambda_{\min}(\Xi^{(n)}),\quad  \mbox{where }\;
\Xi^{(n)}:=\left[ \seq{x_{i}^{(n)},x_{j}^{(n)}} \right]_{i,j=0}^{k_n}\in{\mathbb
    C}^{k_n+1,k_n+1},\quad n>0.
$$
If $\xi>0$ then neither $p_n(\lambda_0)$ nor $\rev p_n(\lambda_0)$ converge to 0 for any $\lambda_0\in\mathbb C$.
\end{proposition}
\proof
For $\lambda_0\in\mathbb C$, with $u^{(n)}=\mat{ccc}  \lambda_0^{0}   \cdots  &   \lambda_0^{k_n} \rix^\top\in{\mathbb C}^{k_n+1}$, we have 
\begin{eqnarray*}
\big\|p_n(\lambda_0)^2\big\|^2&=& \sum_{i,j=0}^{k_n} \overline{\lambda_0}^{i}\lambda_0^{j}\seq{x_{i}^{(n)},x_{j}^{(n)}} \\
&= & u^{(n)*} \Xi^{(n)}  u^{(n)}\\
&\geq & \big\|u^{(n)}\big\|^2 \lambda_{\min}(\Xi^{(n)})\\
&\geq & \xi,\quad n\geq 0.
\end{eqnarray*}
Similarly, we obtain
\begin{eqnarray*}
\big\|\rev p_n(\lambda_0) \big\|^2&=& 
u^{(n)*} \left[ \seq{x_{k_n-i}^{(n)},x_{k_n-j}^{(n)}}\right]_{ij=0}^{k_n} u^{(n)}\\
&\geq & \big\|u^{(n)}\big\|^2 \lambda_{\min}(S\,\Xi^{(n)}S)\\
&= & \big\|u^{(n)}\big\|^2 \lambda_{\min}(\Xi^{(n)})\\
&\geq & \xi,\quad n\geq 0,
\end{eqnarray*}
where $S$ denotes the permutation matrix that reverses the order of indices. 
\eproof

Note that the condition provided  in Proposition~\ref{PP} is sufficient, but not necessary.

\begin{example}\rm
The (constant) sequence of $\mathbb C^3$-valued polynomials 
\[
p_n(\lambda)=p(\lambda)=e_1+(\lambda+\lambda^2) e_2+\lambda^3 e_3
\]
clearly satisfies
$\lim_{n\to \infty} p_n(\lambda_0)\neq 0$ and $\lim_{n\to \infty}\rev p_n(\lambda_0)\neq 0$ for all $\lambda_0\in\mathbb C$, but for all $n>0$ the matrix 
$\Xi^{(n)}$ is singular.
\end{example}

\begin{example}\label{ex:approxchain}\rm 
Let $\mathcal X=\bigoplus_{n=1}^\infty\mathcal X_n$, where $\mathcal X_n=\mathbb C^{2n+1,2n+1}$
and define the operator pencil $P(\lambda)=\lambda E-A$ with
\[
E=\bigoplus_{n=1}^\infty E_n,\; A=\bigoplus_{n=1}^\infty A_n,\quad E_n=\mat{ccc}0&0&I_n\\ 0&\alpha_n&0\\ I_n&0&0\rix,\ A_n=\mat{ccc}\alpha_n &0&0\\ 0&0&I_n\\ 0&I_n&0\rix,
\]
where $\alpha_n$ ($n\geq 1$) is a family of nonnegative parameters. First note that if $\alpha_n=0$ for some $n$ then the pencil is point singular, while if $\inf_{n\geq 1}\alpha_n>0$ then the pencil is not approximately singular due to Lemma~\ref{xns}.
Hence, the interesting case to consider is when $\lim_{n\to \infty}\alpha_n= 0$  but $\alpha_n\neq 0$ for all $n$.

Let $e_j^{(i)}$ denote the vector from $\mathcal X$ that has zero components in all $\mathcal X_k$
different from $i$ and corresponds to the $j$-th standard basis vector of $\mathbb C^{2i+1}$ in
the component $\mathcal X_i$. Then 
\[
p_n(\lambda):=\sum_{j=0}^{n} \lambda^j  e_{j+1}^{(n)},\quad n\geq 1
\]
is a sequence of polynomials, satisfying the assumptions of Proposition~\ref{PP}. Indeed, we have $\Xi^{(n)}=I_n$ and hence $\xi=1$. Furthermore, we have
\[
Ae_1^{(n)}=\alpha_n e_1^{(n)},\quad Ae_{j+1}^{(n)}=e_{n+j+1}^{(n)}=Ee_{j}^{(n)},\;j=1,\dots,n,\quad
Ee_{n+1}^{(n)}=\alpha_n e_{n+1}^{(n)}.
\]
Hence, for arbitrary $\lambda_0\in\mathbb C$
we have
\begin{eqnarray*}
    \norm{(  {\lambda_0} E-A)p_n(\lambda_0)} &=& \|Ae_1^{(n)}\|+\sum_{j=1}^{n}|\lambda_0|^j\|Ee_{j+1}^{(n)}-Ae_j^{(n)}\|+|\lambda_0|^{n+1}\|Ee_{n+1}^{(n)}\|\\
   & =&\alpha_n (1+ |\lambda_0|^{n +1}) 
    =  \norm{(  {\lambda_0} A-E)\rev p_n(\lambda_0)}.
\end{eqnarray*}
Hence, for $\alpha_n=\frac1{(n+1)!}$ the pencil $\lambda E-A$ is approximately singular.
\end{example}

Example~\ref{ex:approxchain} provides, in particular, an instance where the
sequence of degrees of any right approximate polynomial sequence 
is not bounded. However, in general, right approximate polynomial sequences 
may contain right approximate polynomial sequences whose individual
polynomials have smaller degrees than the original ones. Compared to the
concept of singular chains this corresponds to the observation that
linearly dependent singular chains contain singular chains of
smaller length. 

\begin{example}\rm
Let $\mathcal X$ be the Hilbert space from Example~\ref{ex:approxchain} and let $P(\lambda)=\lambda E-A$ with
\[
E=\bigoplus_{n=1}^\infty E_n,\; A=\bigoplus_{n=1}^\infty A_n,\quad E_n=\frac{1}{n}\mat{ccc}1&0&0\\ 0&0&I_n\\ 0&I_n&0\rix,\;A_n=\frac{1}{n}\mat{ccc}0&0&I_n\\ 0&1&0\\ I_n&0&0\rix.
\]
Then as in Example~\ref{ex:approxchain} the polynomials
$$
p_n(\lambda):=\sum_{j=0}^{n} \lambda^j  e_{j+1}^{(n)},\quad n\geq 1
$$
form a right approximate polynomial sequence for $\lambda E-A$.
We also have $Ee_1^{(n)},Ae_1^{(n)}\to 0$ for $n\to\infty$, so also the constant polynomials $e_1^{(n)}$, $n\geq 1$, form a right 
approximate polynomial sequence for $\lambda E-A$.
\end{example}

\section{Operator differential-algebraic equations. Nonuniqueness of solutions. }\label{sec:DAE}

In this section, we consider the effect of singularity of operator pencils on the solutions of the Cauchy problem associated with the
linear operator DAE \eqref{odae}. 
\begin{equation}\label{Cauchy0}
E\dot x(t) = Ax(t), \quad 
x(0) = x_0,
\end{equation}
where  $E\in\mathbf B (\mathcal X,\mathcal Y)$ and    $A$ is a closed and densely defined operator from $\D A\subseteq \mathcal X$ to $\mathcal Y$. Here and everywhere below  $t$ is a real variable and $\dot x$ denotes the partial derivative with respect to $t$.
 We present the notions of  classical solution and mild solutions
to \eqref{Cauchy0}.
\begin{definition}
    \rm A continuous function $x : [0,T) \to \mathcal X$ is called
 a \emph{classical solution} of \eqref{Cauchy0}, if $x$ is continuously differentiable on $[0,T]$, 
$x(t) \in \D A$  for each $t \in[0,T]$  and \eqref{Cauchy0} holds.
Furthermore,  $x(t)$ is called a \emph{mild solution} of~\eqref{Cauchy0}, if it is continuous, $x(0) = x_0$ and for all $ t \in[0,T]$ we have 
$\int_0^t x(s) ds \in \D A$, 
and
$Ex(t) -A\int_0^t x(s) ds
 = Ex_0.
 $
 \end{definition}
    Obviously, a classical solution of   {\eqref{Cauchy0}} is also a mild solution of    {\eqref{Cauchy0}}. We refer the reader to \cite{Tro18} for a wide discussion on consistent initial values and solvability criteria. 
However, note that the results of \cite{Tro18}, as well as many other reults on the subject, require nonemptyness of the set of regular points. In the following we drop this assumption to study singularity of the operator pencil and its relation with 
(non)uniqueness of the solutions of \eqref{Cauchy0}. 

Recall that if $E,A$ are matrices, then the analysis is well understood.  
If the pencil $\lambda E-A$ is regular, then equation \eqref{Cauchy0} has a solution for any consistent initial value, see Theorem 2.12 of \cite{KunM06}. The situation is clearly much more complicated in the operator case, where only sufficient conditions for solvability are known, see e.g. \cite{MehZ23,Tro18} for the case when the resolvent set is nonempty.
We will {\em not} discuss this topic here. 
Another classical fact in the finite dimensional case is the non-existence of solutions of the inhomogeneous DAE $E\dot x =Ax+g$ due to the presence of left singular blocks, see Theorem 2.14(ii) of   \cite{KunM06}.  

It follows from Theorem 2.14 in \cite{KunM06} that  the matrix pencil $\lambda E-A$ has a right singular polynomial if and only if the equation \eqref{Cauchy0} with the initial condition $x(0)=0$ has a nonzero solution. (Observe that existence of such solutions is equivalent to nonuniqueness of solutions of \eqref{Cauchy0}  with nonzero initial condition $x(0)=x_0$, if only such solutions exist.)

We will discuss an infinite dimensional analogue of this fact in the context of the  singularity concepts discussed so far in the paper. We begin with a general result, showing non-uniqueness of classical solutions (hence, also mild solutions) under the assumption of existence of a  generalized singular chain. 
 \begin{theorem}\label{nonunique}
Consider an operator pencil of the form \eqref{Pgeneral} and the associated Cauchy problem~\eqref{Cauchy0}.  Let $\{a_k\}_{k=1}^\infty$ with $a_k\in\D A$ $(k\ge1)$ be a nonzero sequence  such that
\begin{equation}\label{EaAa1}
  Ea_1=0,\quad  Ea_{k+1}=Aa_k,\quad k\geq 1 ,
\end{equation}
and that $\norm{a_k}$, $\norm{Aa_k}$, $\norm{Ea_k}$ are bounded by $\left( \frac kc\right)^k$  for $k$ sufficiently large and some $c>0$. Then there exists a nonzero 
 $\mathcal X$-valued analytic function $f(t)$, defined for $ |t|<1/(ce)$ and satisfying 
\[
f(t)\in\D A,\quad E\dot f(t)=Af(t),\quad f(0)=0.
\]
\end{theorem}
\proof
Define $f(t)$ by the following series
\begin{equation}\label{fcon}
f(t)=\sum_{j=1}^\infty \frac{a_j}{j!}t^j.
\end{equation}
We then observe that the series $\sum_{j=1}^\infty \frac{Xa_j}{j!}t^j$ ($X\in\{I_\mathcal X,A\}$) and $\sum_{j=1}^\infty \frac{Xa_j}{(j-1)!}t^{j-1}$  ($X\in\{I_\mathcal X,E\}$)  are convergent for  $|t|<1/ce$ due to
\[
\limsup_{k\to\infty} \left( \frac{ \left( \frac kc\right)^k } {k!} \right)^{1/k} = \frac1{ce}>0.
\]
Hence, 
$f(t)$ is well defined, and by closedness of $A$ we have $f(t)\in\D A$ for small $t$. 
By definition $f(0)$ is zero and a direct calculation shows $E\dot f(t)=Af(t)$ everywhere where the series defining $f(t)$ converges.
\eproof
Based on this we obtain the relation between the singularity notions for operator  pencils and (non)uniqueness of the solutions of the corresponding Cauchy problem.
\begin{corollary}
    \label{polyL}
    Consider the Cauchy problem  \eqref{Cauchy0}.  If  the operator pencil $\lambda E-A$ has a right singular polynomial,
then there exists a nonzero $\mathcal X$-valued polynomial   {$f(\lambda)$} 
with $f(0)=0$ satisfying 
$f(t)\in\D A$, and $E\dot f(t)=Af(t)$ for $t\in\mathbb R$.
\end{corollary}
\proof
Let $p(\lambda)=a_{k+1}+\cdots+\lambda^ka_1$ be a right singular polynomial of $\lambda E-A$. We show that  after appending the coefficients  with an infinite number of zeros, the sequence $\{ a_k \}_{k=1}^\infty$ satisfies the assumptions of Theorem ~\ref{nonunique}.
  By Lemma~\ref{division} we may assume that $a_{k+1}\neq 0$, hence the sequence  $\{ a_k \}_{k=1}^\infty$ is nonzero. Now we show that $a_k\in \mathcal D(A)$ for $k\geq 1$. 
By definition, for every $\lambda_0\in\mathbb C$ the vector $p(\lambda_0)$ is in $\D{A}$. Taking distinct $\lambda_0,\dots,\lambda_k\in\mathbb C$ we obtain $k+1$ vectors
\[
b_j=\sum_{i=0}^k \lambda_j^i a_i \in\D A,\quad j=0,\dots k.
\]
Consider  the Vandermonde matrix   $\Lambda:=\mat{c} \lambda_j^i\rix_{i,j=0}^{k+1}$ and its inverse  $\Lambda^{-1} =:\mat{c} \mu_{ij}\rix_{i,j=0}^{k+1}$.
Then 
\[
a_j=\sum_{i=0}^k \mu_{ji} b_i \in\D A,\quad j=0,\dots k.
\]
Since $p(\lambda)$ is a right singular polynomial we have, by definition, $(\lambda E-A)p(\lambda)=0$, which immediately gives \eqref{EaAa1}. 
Now observe that \eqref{fcon} is also a polynomial, hence $E\dot p(t)=Ap(t)$ for all $t\in\Real$.   {Choosing now $f(\lambda)=\lambda p(\lambda)$ finishes the proof.}
\eproof

\begin{example}\label{facfac}\rm
  Let $E$ be the closure of the forward shift operator and let $A$ be the closure of a diagonal operator in $\ell^2$ defined by
\[
Ee_1=0, \ Ee_{k+1}=e_k,\quad Ae_k=(k+1)e_k,\quad  (k\geq 1),
\]
extending them as usually to linear operators on the linear span of basis vectors and taking the closure. Then $E$ is bounded but $A$ is not.
  The sequence  $a_k=k! e_k$ satisfies the assumptions of Theorem~\ref{nonunique}. On the other hand  $A$ is invertible and hence $\lambda E-A$ neither has left nor right singular polynomials.
\end{example} 

We also have a relation between non-uniqueness of solutions of \eqref{Cauchy0} and  the existence of a right singular function for $\lambda E-A$.
\begin{remark}
  {\rm  If $\norm{a_k}$, $\norm{Aa_k}$, $\norm{Ea_k}$ are bounded by  $(k/c)^{-k}$  for $k$ sufficiently large and some $c>0$ then $x(\lambda)=\sum_{j=1}^\infty a_j\lambda^{-j}$ is a right singular function of  $\lambda E-A$,  as in Section~\ref{ss:singularanalytic}. However, for $A,E$ as in Example~\ref{facfac} the series $x(\lambda)=\sum_{j=1}^\infty a_j\lambda^{-j}$ diverges for all $\lambda\neq 0$,  }
\end{remark}

Finally, we discuss whether non-uniqueness of the solutions of \eqref{Cauchy0} is connected with   approximate singularity of $\lambda E-A$. Apparently, there does not seem to be any relation, see the following two examples.

\begin{example}\label{ex:36}\rm
Let $E$ be the shift on the canonical basis $e_1,e_2,\dots$ of $\ell^2$, i.e., $Ee_k=e_{k-1}$,$Ee_1=0$ and let $A=I_{\ell^2}$.
Then $E$, $A$ and $a_k=e_k$ ($k\geq 1)$ satisfy the setting of Proposition~\ref{nonunique} and $f(t)=\sum_{j=1}^\infty \frac{e_j}{j!}t^j$ is a solution of $E\dot f=Af$, with $f(0)=0$. 
On the other hand $0$ is clearly a regular point of the 
pencil $\lambda E-A$. In particular, the pencil is not approximately singular. 

{  Furthermore, note that there exist nonzero initial values $g(0)$ for which the equation $Eg'=Ag$ has a solution, which is then necessarily not unique. For example,  the initial condition  $g(0)= e_1$ produces solutions of the form $g(t)=e_1+\sum_{j=2}^\infty\frac1{(k-1)!}e_kt^{k-1} +\alpha f(t)$, where $\alpha\in\Real$ and $f$ is as above}.
\end{example}

\begin{example}\label{ex:37}\rm
Consider $A=E=\diag(1,1/2,1/3,...)$ in $\ell^2$ from Example~\ref{ex:diag}. The DAE $E\dot x= Ax$ has precisely one solution for any initial condition $x(0)=x_0\in\ell^2$, namely $x(t)=  {e^{t}}x_0$.
\end{example}

In this section we have studied the relationship between the different singularity concepts and the uniqueness of solutions of the Cauchy problem~\eqref{Cauchy0}. In the next section we discuss the particular class of dissipative Hamiltonian operator pencils.

\section{Dissipative Hamiltonian Operator pencils}\label{sec:dhpencils}

As discussed in the introductory section, a major motivation for the analysis of operator DAEs is the
study of dissipative Hamiltonian operator pencils. We will consider these pencils in the slightly more general form
\begin{equation}\label{PEJR}
P(\lambda)=\lambda E - BQ,
\end{equation}
where 
\begin{enumerate}[\rm (i)]
    \item $E\in{\mathbf B}(\mathcal X,\mathcal Y)$, $Q\in\mathbf B(\mathcal X,\mathcal Y)$ is   {boundedly invertible} and  $Q^*E$ is selfadjoint and nonnegative, i.e., $\seq{Q^*Ex,x}\geq 0$ for all $x\in\mathcal X$;
    \item  $B$ is a closed, densely defined and dissipative operator in $\mathcal Y$, i.e., $\RE\seq{Bx,x}\leq0$ for all $x\in \D B\subseteq \mathcal Y$.
\end{enumerate} 
 Such pencils satisfy the assumptions introduced in Section~\ref{sec:prelim}, as $BQ$ is closed and densely defined.
 In particular, the set of regular points is open (possibly empty).
We note that this class is more general than the class studied  in \cite{MehZ23}, as we neither  assume  maximal dissipativity of $B$ nor that $\ker E$ needs to be an invariant subspace for $Q$.    

\subsection{Singularity of dissipative Hamiltonian operator pencils}\label{sec:singdhadae}
In \cite{MehMW21} singularity of dissipative Hamiltonian pencils was investigated
in the finite-dimensional case.
In particular, it was shown that the presence of eigenvalues in 
the open right half plane of matrix pencils of the form~\eqref{PEJR} already
implies singularity of the pencil, which furthermore
is equivalent to all three matrices $E$, $J$ and $R$ having a common kernel.

In this section we extend these results to the infinite dimensional case. By doing so we obtain an essentially simplified framework for studying singularity,  compared to Section~\ref{sec:sinpen}.

In a recent paper \cite{MehZ23} regularity and singularity of related pencils is studied, though in a slighty different operator setting. The results of \cite{MehZ23} are, however, similar to  Theorems~\ref{tEJR}, \ref{tEJR2} and \ref{tEJR3} below: regularity is equivalent to the open right half-plane being
contained in the set of regular points. 
\begin{proposition}\label{prop14}
Consider an operator pencil of the form~\eqref{PEJR},  let $\lambda_0\in\mathbb C$ 
with $\RE\lambda_0>0$ and
let vectors $x_n\in \mathcal{D}(BQ)$, $n\in\mathbb N$, form a bounded 
sequence. Then for  $n\to\infty$, we have
\begin{equation*}
\label{===}
 P(\lambda_0)x_n\to 0\quad \mbox{ if and only if }\quad  Ex_n\to 0 \mbox{ and } BQ x_n \to 0.
\end{equation*}
In particular, we have
\begin{equation}
\label{==}
\ker P(\lambda_0)=\ker E\cap\ker(BQ) .
\end{equation}

If, additionally, $\mathcal X=\mathcal Y$,  $Q=I_{\mathcal X}$ and $B=J-R$, with both $J$ and $R$ being closed, $\D J\subseteq \D R$ or $\D R\subseteq \D J$ and $\seq{Jx,x}\in\ii\Real$, $\seq{Rx,x}\geq 0$ for $x\in\D R\cap \D J$ then 
\begin{equation}
\label{==1}
\ker P(\lambda_0)=\ker E\cap\ker J\cap \ker R=\ker(E^2+R^2-J^2).
\end{equation}
\end{proposition}
\proof Let $ P(\lambda_0)x_n\to0$, i.e.,  $(\lambda_0 E-BQ)x_n\to0$. Then, $(\lambda_0 Q^*E-Q^*BQ)x_n\to0$
and  $\seq{(\lambda_0 Q^*E-Q^*BQ)x_n,x_n}\to 0$. Taking the real part  and using that
$\seq{\big(\RE(\lambda_0) Q^*E\big)x_n,x_n}\geq 0$ as well as $\RE\seq{(Q^*BQ)x_n,x_n}=\RE\seq{BQx_n,Qx_n}\leq 0$, we obtain
$\seq{Q^*Ex_n,x_n}\to0$.  The Cauchy-Schwarz
inequality for the semidefinite inner product $\seq{Q^*E\cdot,\cdot}$ (cf. \cite[Theorem 4.2]{Ru87}) and the boundedness of $Q^*E$ then yield that
\begin{equation}\label{E2E}
\seq{Q^*Ex_n,Q^*Ex_n}\leq \seq{Q^*Ex_n,x_n}^{1/2}\seq{(Q^*E)^2x_n,Q^*Ex_n}^{1/2}\leq \seq{Q^*Ex_n,x_n}^{1/2}\cdot\norm {Q^*E}^{3/2}\cdot\|x_n\| .
\end{equation}
Hence, $Q^*Ex_n\to 0$ and thus $Ex_n\to0$, since $Q$ is   {boundedly invertible}, and $BQx_n\to 0$. The converse implication is trivial and to see \eqref{==} take $x\in \ker P(\lambda_0)$ and a constant sequence $x_n=x$. 

To prove the second part of the claim, let $x\in\ker P(\lambda_0)$. Due to the first part  we
 obtain $Ex=0$ and $(J-R)x=0$.  Since $\seq{Jx,x}$ is purely imaginary and
$\seq{Rx,x}$ is real, we have $\seq{Rx,x}=0$. Recall that  $R$ is positive semidefninite (not necessarily bounded), hence, $|\seq{Rx,y}|\leq \seq{Rx,x}^{1/2}\seq{Ry,y}^{1/2}=0$ for any $y\in\D R$. Consequently $x\in\ker R$
 and finally $Jx=0$ as well, which shows that $x\in\ker E\cap\ker R\cap\ker J$. The
 inclusion $ \ker E\cap\ker R\cap\ker J\subseteq \ker(E^2-J^2+R^2)$ is trivial.  If now $x\in\ker(E^2-J^2+R^2)$ then $\seq{(E^2-J^2+R^2)x,x}=0$ and due to the fact that $E^2,R^2$ and
  $-J^2$ are positive semidefinite we have $Ex=Rx=Jx=0$, which shows $x\in\ker P(\lambda_0)$. This finishes the proof of
  \eqref{==1}.
\eproof

\begin{remark}\label{rem:Q=I}{\rm
If we are in the special case  that $Q=I_{\mathcal X}$ in Proposition~\ref{prop14}  and if the coefficient $R$ of the pencil $P(\lambda)$  
is not only symmetric, but also bounded (and hence self-adjoint),
then the condition $QBx_n=(J-R)x_n\to 0$ is in fact equivalent to $Jx_n,Rx_n\to 0$.
Indeed, $(J-R)x_n\to 0$ implies $\langle Jx_n,x_n\rangle$, $\langle Rx_n,x_n\rangle\to 0$ and from the latter we obtain $Rx_n\to 0$ by the same reasoning as in~\eqref{E2E}.
But then we must also have $Jx_n\to 0$. Without the boundedness of $R$ one can  show weak convergenve of $Rx_n$ to $0$, i.e., $\seq{Rx_n,y}\to0$ for all $y\in\mathcal X$.
}
\end{remark}

The next results investigate the effects of the presence of singular points in the
open right half plane. We will distinguish between point singularities,
approximate singularities and general singular points. Let us point out the essential difference with the general case discussed in Section~\ref{sec:sinpen}. From Theorem~\ref{tEJR} it follows that $s_{p}\big(P(\lambda)\big)=\mathbb C$ \emph{is} equivalent to being point singular (cf.~in contrast Example~\ref{ex:nodecom}), and from Theorem~\ref{tEJR2} it follows that $s_{ap}\big(P(\lambda)\big)=\mathbb C$ \emph{is} equivalent to being approximately singular,  cf. Example~\ref{ex:nondecom2}). 

\begin{theorem}\label{tEJR}
Consider an operator pencil $P(\lambda)$ of the form~\eqref{PEJR}. Then the following 
conditions for point singularities are equivalent:
\begin{enumerate}[\rm (a)]
\item $s_{\p} \big(P(\lambda)\big)=\mathbb C$ and $\infty$ is a point singularity;
\item\label{s:IM} $s_{\p} \big(P(\lambda)\big)\cap\set{\lambda_0\in{\mathbb C}:\RE\lambda_0>0}
    \neq\emptyset$;
\item $\ker E\cap\ker(BQ)\neq\set0$;
\item the pencil $P(\lambda)$ has a right singular polynomial.
\end{enumerate}
\end{theorem}
\proof
The implication (a) $\To $ (b)  is trivial, (b) $\To$ (c) follows from Proposition~\ref{prop14}. To see (c) $\To$ (d)  let $p(\lambda)=x$, where $x$ is a nonzero element of  $\ker E\cap\ker BQ$. Clearly it is a right singular polynomial, see Definition~\ref{def:PS}.   The implication (d) $\To$ (a) is a special instance of Proposition~\ref{polyC}.
\eproof

For approximate singularities we have the following equivalence result.
\begin{theorem}\label{tEJR2}
Consider a pencil $P(\lambda)$ of the form \eqref{PEJR}. Then the following 
conditions for approximate singularities are equivalent:
\begin{enumerate}[\rm (a)]
\item $s_{\ap} \big(P(\lambda)\big)=\mathbb C$  and $\infty$ is an approximate singularity;
\item $s_{\ap}\big( P(\lambda)\big)\cap\set{\lambda_0\in{\mathbb C}:\RE\lambda_0>0}
    \neq\emptyset$;\label{s:IM2}
\item there exists a sequence of unit norm vectors $x_n\in \mathcal D (BQ)$ 
such that $Ex_n,BQx_n\to 0$;
\item the pencil $P(\lambda)$ has a right  approximate singular polynomial sequence.
\end{enumerate}
\end{theorem}
\proof
The implication (a) $\To$ (b) is trivial,  (b) $\To$ (c) follows from Proposition~\ref{prop14}.  The implication (c) $\To$ (d) follows from the definition and (d) $\To$ (a) is a special case of Theorem~\ref{PP?}.
\eproof

{   Finally, we present a result on the set of singular points, cf. \cite{MehZ23} for a similar result under stronger assumptions. For this purpose we need a notion of a \emph{maximally dissipative} operator, that is, 
a dissipative operator $B$ such that  its  graph is not contained in the graph of some other  dissipative operator. This is equivalent to saying that $B$ is dissipative and the range of $B$ is the whole space, and also equivalent to saying that $B$ and $B^*$ are both dissipative, see, e.g., the Appendix of \cite{MehZ23} for a convenient summary. 

\begin{theorem}\label{tEJR3}
 Consider a pencil $P(\lambda)$ of the form \eqref{PEJR}. Then the following conditions are equivalent:
\begin{enumerate}[\rm (a)]
\item  $s \big(P(\lambda)\big)\neq \mathbb C$ and $B$ is maximally dissipative;
\item $s\big( P(\lambda)\big)\cap\set{\lambda_0\in{\mathbb C}:\RE\lambda_0>0} =\emptyset$;\label{s:IM3}
\item $\set{\lambda_0\in{\mathbb C}:\RE\lambda_0>0} \setminus s\big( P(\lambda)\big) \neq\emptyset$.
\end{enumerate}
\end{theorem}
\proof
(a)$\To$(b) Assume (a), and suppose that $\mu_0\in s\big(P(\lambda)\big)$, $\RE\mu_0>0$. Point (a)  implies  that the set of approximate singularities is not the whole complex plane.
By Theorem~\ref{tEJR2}, we have that 
 $s_{\ap}\big( P(\lambda)\big)\cap\set{\lambda_0\in{\mathbb C}:\RE\lambda_0>0}=\emptyset$, in particular $\mu_0\notin s_{\ap}(P(\lambda))$. 
In other words, the range of $\mu_0 E- BQ$ is closed, 
but not equal to the whole space. Therefore, by the kernel-range decomposition,
\begin{equation}\label{ker*}
\ker \big(\bar\mu_0 E^* - (BQ)^*\big)\neq\{0\}.
\end{equation}
As $B$ is maximally dissipative,  $B^*$, and thus also $Q^*B^*Q$, is dissipative. Therefore, the  operator pencil
$$
\widetilde P(\lambda)=\lambda E^* - (BQ)^*= \lambda E^* -Q^*B^* =  \lambda E^* -(Q^*B^*Q) Q^{-1},
$$
satisfies \eqref{PEJR} with assumptions (i), (ii), where $E^*$ plays now the role of $E$, $Q^*B^*Q$ of $B$ and $Q^{-1}$ of $Q$. By Theorem~\ref{tEJR} and \eqref{ker*} we have that $s_{\p}(\widetilde P\big(\lambda)\big)=\mathbb C$ and thus, again by therange-kernel argument, $s\big(P(\lambda)\big)=\mathbb C$, contradiction.

 The implication (b)$\To$(c) is trivial, we show now (c)$\To$(a). Assume that $\lambda_0 E- BQ$ is   {boundedly invertible} for some $\lambda_0$ with $\RE\lambda_0>0$. Then $-\lambda_0 Q^*E + Q^*BQ$ is   {boundedly invertible} as well, furthermore it is dissipative.
 Hence, it is maximally dissipative, in consequence $Q^*BQ$, and thus $B$, is maximally dissipative.
\eproof

}

\begin{remark}\label{rem:notequi}{\rm
We have shown that conditions \eqref{s:spC}, \eqref{s:sapC}, and \eqref{s:sC}  are not equivalent for general operator pencils. In fact, they are also not equivalent for operator pencils of the form~\eqref{PEJR}. Suitable counterexamples are essentially the same ones as before. For the first situation take any pencil $\lambda E-E$ with $E$ selfadjoint bounded, positive
semidefinite, with zero in the approximate spectrum but not in the point spectrum. For the second one
take the pencil in \eqref{tildeA} with $\tilde A=\ii B$, where $B$ is symmetric but not
selfadjoint.
}
\end{remark}

Using Theorem~\ref{distto0} and Theorem~\ref{tEJR2} we have  the following corollary.
\begin{corollary}
If $\lambda E-BQ$ is a pencil as in \eqref{PEJR} satisfying one of the equivalent
conditions in Theorem~\ref{tEJR2} then for any sequence of projections
$P_n$, $n\geq1$ satisfying \eqref{P1}--\eqref{P3} the distance to singularity of the finite dimensional operator
pencils
$
P_n ( \lambda E- A  ) P_n
$
converges to zero for $n\to\infty$.
\end{corollary}

\subsection{Uniqueness of solutions of the corresponding operator DAEs}
Consider the dissipative operator DAE
\begin{equation}\label{Cauchy1}
E\dot x(t)=BQx(t),\quad x(0)=x_0, \quad t\in[0,T],
\end{equation}
with $E,B,Q$ as in~\eqref{PEJR}. We refer to Section~\ref{sec:DAE} for the definition of classical and mild solutions. 
It was also presented there that neither nonemptyness of the set of regular points of $\lambda E-A$ implies uniqueness of the solutions of \eqref{Cauchy1} nor conversely, see Example~\ref{ex:36} and 
Example~\ref{ex:37}, respectively. While Example~\ref{ex:37} is already of the dissipative Hamiltonian form~\eqref{PEJR}, here $Q=I_{\ell^2}$, Example~\ref{ex:36} cannot be written in this form. This will become apparent, as we will show that  for dissipative Hamiltonian pencils with the set of point singularities not equal to the whole complex plane the solution of \eqref{Cauchy1} is unique (if it exists). This is a typical fact in energy based approaches, see, e.g., ~\cite{Eva22}. The key result that is needed for this is the power balance equation, see \cite{Mor24} for a detailed study,
\begin{equation}\label{PBE}
    \seq{Ef(t_0),Qf(t_0)} - \seq{Ef(0),Qf(0)} = 2\RE\int_0^{t_0} \seq{Q^*BQ f(t),f(t)} dt,\quad { t_0\in[0,T]},
\end{equation}
which holds for any classical solution of~\eqref{Cauchy1}. Indeed, if $f(t)$ is a classical solution, then
\[
\frac d{dt}\seq{Ef(t),Qf(t)}=2\RE\seq{E\dot f(t),Qf(t)}=2\RE\seq{Q^*BQ f(t),f(t)},
\]
which after integration gives~\eqref{PBE}.
Note that Theorem~\ref{th:ce} provides in fact another equivalent condition to the conditions (a)--(d)
from Theorem~\ref{tEJR}. However, for presentation reasons, we work with a negated version of (c).

\begin{theorem}\label{th:ce}
    Consider a linear operator pencil of the form~\eqref{PEJR}. Then the following statements are equivalent:
    \begin{enumerate}
    \item[\rm (c')] $\ker E\cap\ker BQ= \{0\}$;
 \item[\rm (e')]
 if $x_1(t),x_2(t)$ are two mild solutions of \eqref{Cauchy1} (with the same initial condition), then $x_1(t)=x_2(t)$ for $t\in[0,T]$. 
    \end{enumerate}
\end{theorem}

\proof To show  (e')$\Rightarrow$(c') assume that (e') holds and that we have $\ker E\cap\ker BQ\neq \{0\}$.
Hence, there exists a nonzero continuously differentiable function $x$ satisfying $x(0)=0$ and $x(t)\in\ker E\cap \ker (BQ)$ for all $t\in[0,T]$. Assume that $x_1(t)$ is a mild solution  of \eqref{Cauchy1}. 
Then $x_2(t)=x_1(t)+x(t)$ is a mild solution of \eqref{Cauchy0} as well. Further $x_2(t)$ and $x_1(t)$ are not  equal almost everywhere, which is a contradiction.

To show the converse implication assume that  $\ker E\cap\ker (BQ)= \{0\}$. Suppose that $x_1(t),x_2(t)$ are two mild solution of \eqref{Cauchy1} with the same initial condition.
Hence,  $u(t)=x_1(t)-x_2(t)$ is a mild solution of \eqref{Cauchy1} with the initial condition $u(0)=0$. Therefore, $f(t)=\int_0^t u(s)ds$ is a classical solution of \eqref{Cauchy1} with the initial condition  $f(0)=0$. Hence, \eqref{PBE} holds, and for any $t_0\in[0,T]$ the left hand side  is nonnegative while the right hand side is nonpositive. Thus, they are both zero. In particular, 
$\seq{Q^*E f(t_0),f(t_0)}=0$, which gives $Ef(t_0)=0$, since $Q$ is   {boundedly invertible}. On the other hand,
$2\RE\int_0^{t_0} \seq{Q^*BQ f(t),f(t)} dt=0$. As $f(t)$ is continuous and $\ker(BQ)=\ker(Q^*BQ)$ is a closed space, we have that $BQf(t_0)=0$. Hence,
$f(t_0)\in\ker E\cap\ker BQ=\{0\}$, by (c'). Consequently,  $x_1(t_0)=x_2(t_0)$ for $t_0\in[0,T]$. 
\eproof

\section{Conclusions}
The spectral theory for operator pencils has been studied. Three different concepts of singular operator pencils have been introduced and analyzed in detail.  Many examples illustrate the subtle differences. The results are then applied to operator pencils arising in dissipative Hamiltonian differential-algebraic equations, and it is shown that many results known from the finite dimensional case of matrix pencils extend directly to the infinite dimensional setting.

\section*{Acknowledgment}

 The authors are grateful to Jan Stochel for valuable ideas and discussions, {  they also express their gratitude to the anonymous referees for their careful reading of the manuscript and for many helpful comments.}

V.M. has been supported by Deutsche Forschungsgemeinschaft (DFG) through the SPP1984 ``Hybrid and multimodal energy systems'' Project: \emph{Distributed Dynamic Security Control} and by the DFG Research Center Math+, Project  No.~390685689. 
\emph{Advanced Modeling, Simulation, and Optimization of Large Scale Multi-Energy Systems}.

 M.W. was supported by   the Jagieloinan University ID.UJ 
grant, decision No. PSP: U1U/P06/NO/02.12

\bibliographystyle{plain}
\bibliography{mmw}

\end{document}